\documentclass[12pt]{article}

\usepackage{t-angles}
\usepackage{amsmath}

\usepackage{amsfonts}
\usepackage{latexsym}
\usepackage{amssymb}

\newcommand{\id}{\operatorname{id}}
\newcommand{\ev}{\operatorname{ev}}

\newcommand{\ot}{\otimes}

\newcommand{\trl}{\triangleleft}
\newcommand{\trr}{\triangleright}

\newcommand{\di}{{}_{(1)}}
\newcommand{\dii}{{}_{(2)}}

\newcommand{\lcom}[1]{{}^{#1}\!{\cal M}}           

\newcommand{\lyd}[1]{{}^{#1}_{#1}{\cal YD(C)}}    
\newcommand{\ryd}[1]{{\cal YD(C)}^{#1}_{#1}}     

\def\rbiprod{{\cdot\kern-.33em\triangleright\!\!\!<}}
\def\lbiprod{{>\!\!\!\triangleleft\kern-.33em\cdot\, }}

\def\lrbiprod{{\ \cdot\kern-.60em\triangleright\kern-.33em\triangleleft\kern-.33em\cdot\, }}

\def\lprod{{>\!\!\!\triangleleft\kern-.33em\ \, }}
\def\rprod{{\triangleright\!\!\!<}}

\newcommand{\lcoprod}{{>\!\!\blacktriangleleft}}
\newcommand{\rcoprod}{{\blacktriangleright\!\!<}}

\newcommand{\lrcoprod}{{\blacktriangleright\!\!\blacktriangleleft}}

\newcommand{\lcprp}{{\vartriangleright\!\!\blacktriangleleft}}
\newcommand{\lprcp}{{\blacktriangleright\!\!\vartriangleleft}}

\title{Double Bicrosssum of Braided Lie algebras}
\author{ Shouchuan Zhang \ \ Tao Zhang\\
Department of Mathematics, Hunan University\\
Changsha 410081,  R.P.China}
\date{}
\date{}
\begin{document}
\newtheorem{Theorem}{\quad Theorem}[section]
\newtheorem{Proposition}[Theorem]{\quad Proposition}
\newtheorem{Definition}[Theorem]{\quad Definition}
\newtheorem{Corollary}[Theorem]{\quad Corollary}
\newtheorem{Lemma}[Theorem]{\quad Lemma}
\newtheorem{Example}[Theorem]{\quad Example}
\newtheorem{Remark}[Theorem]{\quad Remark}

\maketitle \addtocounter{section}{-1}

\begin {abstract}
The condition for  double bicrosssum  to be a braided Lie bialgebra
is given. The result generalizes quantum double, bicrosssum,
bicrosscosum, bisum. The quantum double of braided Lie bialgebras is
constructed. The relation between double crosssum of  Lie algebras
and double crossproduct of Hopf algebras is given.

\end {abstract}



\section{Introduction}

Lie bialgebras and quasitriangular Lie bialgebras were introduced by
Drinfeld in his remarkable article \cite{Dr87}. In the sequel,  many
important examples of Lie bialgebras were found. For example, Witt
algebra and Virasoro algebra in \cite{Mi94} and \cite{Ta93};
parabolic Lie bialgebra in \cite{Ma00}; quaternionic Lie bialgebra
in  \cite{Fe05}, and so on. Low dimensional coboundary
 (or triangular ) Lie bialgebras
are classified  in \cite{Zh98}.  Lie bialgebras  play an important
role in quantum groups and Lie groups ( see \cite{Ma95, CP94}).
Majid studied the bisum of Lie bialgebras in \cite{Ma90b} and
\cite{Ma00}. Grabowski studied the double bosonisation of Lie
bialgebras \cite{Gr05}.

Braided tensor categories become more and more important. In
mathematical physics, Drinfeld and Jambo found the solutions of
Yang-Baxter equations by means of braided tensor categories ; In
string theory, Yi-Zhi Huang and Liang Kong  gave the relation
between Open-string vertex algebras and braided tensor categories
(see \cite {Hu05, HK04}); In topology, there exits a tight contact
between braided tensor categories and the  invariants of links and
3-dimensional manifolds ( see \cite {BK01, He91,Ka97, Ra94} ).

In this paper,  we follows   Majid's  work in papers \cite{Ma90b,
Ma00} and generalize  his conclusion in two ways. On the one hand,
we give the condition for double bicrosssum to become  a Lie
bialgebra, which generalizes quantum double, bicrosssum,
bicrosscosum and bisum. We build the relation between the bisum of
Lie algebras and the biproduct of Hopf algebras. On the other hand,
our conclusions hold in braided tensor categories, so our work will
supply a useful tool for conformal fields, vector algebras and
3-manifolds.

\section*{\bf Preliminaries}

We begin with recalling additive category ( see \cite [Page 148]
{Fa73}).

\begin {Definition} \label {1.1'}
 A category  ${\cal C}$ is called an additive category, if
  the following conditions are satisfied for any 如果对 中的任何三个目标
  $U, V, W \in {\cal C}$:

 (Ad1) $Hom _{\cal C}(U,V)$ is an additive group and
$Hom _{\cal C}(V,U) \times Hom _{\cal C}(W, V) \rightarrow Hom
_{\cal C} (W, U)$
 is bilinear;

 (Ad2) ${\cal C} $ has a zero object 0;

 (Ad3) There exists the direct sum for any finite objects in ${\cal
 C}$.
\end {Definition}

\begin {Example} \label {1.1''}
(i) If  $H$ is a Hopf algebra with an invertible antipode, then the
Yetter-Drinfeld category $^H_H {\cal YD} $ is an additive category;

(ii) The braided tensor category, which is studied  in \cite {Hu05},
is an  additive category.
\end {Example}
Otherwise, the most of important examples in \cite {BK01} are
additive braided tensor categories.

In this paper,  unless otherwise  stated, $({\cal C}, C)$ is an
additive braided tensor category with braiding $C$.

 If ${\cal E} $ is a family of objects in  ${\cal C}$ and
 $C_{U, V} = (C_{V, U})^{-1}$  for any $U, V \in {\cal
 E}$, then the braiding $C$ is said to be symmetric on ${\cal E}$.
 If  $C$  is symmetric on $L$ and $L$ has a left duality $L^*$, then, by \cite {Zh03},
 $C$  is symmetric on the set $\{L, L^*\}$.

Let $L$, $U$, $V$, $A$ and $H$  be objects in ${\cal C}$. $C$
 is symmetric on  $L$, $U$and $V$, respectively;  $C$ is symmetric on set
  $\{A,  H\}$;
$\alpha : \ H \otimes A \rightarrow A,
$$ \ \ \beta : \ H \otimes A \rightarrow H,
$$ \ \ \phi \ \ \ A \rightarrow H \otimes A$\ and  \ $\psi \ H
\rightarrow H \otimes A$ are morphisms in ${\cal C}$. We work on the
five objects above throughout this paper.

\begin {Definition} \label {1.1}
  (\cite{Ma88} )\quad $(L, [\ , \ ])$ is called a braided Lie algebra in $\cal
C$, if  $[\ , \ ]: L\ot L\to L$ is a morphism in ${\cal C}$, and the
following conditions are satisfied:

(L1)\ ${\cal C}$-anti- symmetry: $[\ ,\ ](\id+C)=0$,

(L2)\ $\cal C$-Jacobi identity: $[,
]([,]\ot\id)(\id+C_{12}C_{23}+C_{23}C_{12})=0$, where $C_{12} =
C_{L,L} \otimes id_L$ and  $C_{23} = id _L \otimes C_{L,L}$.

\end {Definition}

If $(\Lambda, \beta)$ is a coquasitriangular Hopf algebra, then the
category $ {\cal C} = \lcom \Lambda$ of all left $\Lambda$-
comodules is a symmetric braided tensor category. $\beta$-Lie
algebras in \cite{BFM96} and \cite{BFM01}, as well as $\epsilon$-Lie
algebras in \cite{Sc79}, are braided Lie algebras in ${\cal C}$.

\begin {Definition} \label {1.2} \quad  $(L, \delta)$ is called a braided Lie coalgebra in
$\cal C$, if  $\delta: L\to L\ot L$是${\cal C}$ is a morphism in
$\cal C$ and the following conditions are satisfied:

(CL1)\ $\cal C$-co-anti-symmetry: $(\id+C)\delta=0$,

(CL2)\  $\cal C$-co-Jacobi identity:
$(\id+C_{12}C_{23}+C_{23}C_{12})(\delta \ot \id)\delta=0 $.

\end {Definition}
In the definition above, (L1), (L2), (CL1) and  (CL2) are  also
called braided anti- symmetry, braided Jacobi identity, braided
anti- co-symmetry, braided co-Jacobi identity.

\begin {Example} \label {4.8} \footnote
 {This is given  by Dr. Lian Kong } Assume that  $U$ is an object in ${\cal C}$ with a left
 duality
  $U^*$ and  $X = U \otimes U^*
$, $Y = U ^*\otimes U$. Define the multiplication and unit in $X$ as
follows: $m_X = id _U \otimes ev_U \otimes id _{U^*}$, $\eta _X =
coev _U$. It is clear that  $(X, m, \eta)$ is an associative algebra
in ${\cal C}$. Therefore,  $X$ becomes a braided Lie algebra under
the following bracket  operation: $[\ \ ] = m_X - m _XC_{X, X}$.
Dually, define the comultiplication and counit in $Y$ as follows:
$\Delta_Y = id _{U^*} \otimes coev_U \otimes id _U$, $\epsilon _Y =
ev _U$. It is clear that  $(Y, \Delta , \epsilon )$ is a braided
coalgebra in ${\cal C}$. Therefore,  $Y$ is a braided Lie coalgebra
under the following operation:  $\delta = \Delta _X - C_{Y, Y}\Delta
_Y$.
\end {Example}

\begin {Definition} \label {1.3} \quad Assume that   $(H, [,])$ is a braided Lie algebra in $\cal
C$ and  $(H, \delta)$ is a braided Lie coalgebra in $\cal C$. If the
following condition  is satisfied: \begin {eqnarray*} (LB): \ \ \ \
\ && \delta
[,]=(([,]\ot \id_H)(\id_H\ot \delta)(\id_H \otimes id_H -C_{H,H})\\
&+&(\id_H\ot [,])(C_{H,H}\ot \id_H)(\id_H\ot \delta))(\id_H\ot
\id_H-C_{H, H})\ ,
\end {eqnarray*}
then $(H, [,], \delta)$ is called a braided Lie bialgebra  in $\cal
C$.

\end {Definition}

In order to denote operations for conveniency, we use braiding
diagrams. Let
\[
 \step[2] [\ \ ]=\enspace
\begin{tangle}
\object{H}\step[2]\object{H}\\
\cu*\\
\step\object{H}\\
\end{tangle}
\enspace, \step[2] \delta=\enspace
\begin{tangle}
\step\object{H}\\
\cd*\\
\object{H}\step[2]\object{H}\\
\end{tangle}
\enspace, \alpha=\enspace
\begin{tangle}
\object{H}\step[2]\object{A}\\
\Lu\\
\step[2]\object{A}\\
\end{tangle}
\enspace, \step[2] \beta=\enspace
\begin{tangle}
\object{A}\step[2]\object{H}\\
\Ru\\
\object{A}\\
\end{tangle}
\enspace, \step[2] \phi=\enspace
\begin{tangle}
\step[2]\object{A}\\
\Ld\\
\object{H}\step[2]\object{A}\\
\end{tangle}
\enspace, \step[2] \psi=\enspace
\begin{tangle}
\object{A}\\
\Rd\\
\object{A}\step[2]\object{H}\\
\end{tangle}
\enspace.
\] In particular, when  $C$  is symmetric on  $\{U, V\}$,  $C _{U, V}$ is
denoted by the following diagram:
 $C_{U,V}=\enspace
\begin{tangle}
\object{U}\step[2]\object{V}\\
\XX\\
\object{V}\step[2]\object{U}\\
\end{tangle}
\enspace.$ When  $U$ has a left duality  $U^*$, denote the
evaluation $ev_U =d_U$ and
 coevaluation $coev_U = b_U$ by the following diagrams: \ \ \ \ $ ev _U=\enspace
\begin{tangle}
\object{U^*}\step[2]\object{U}\\
\ev
\end{tangle}
\enspace, coev _U =\enspace
\begin{tangle}
\coev \\
\object{U}\step[2]\object{U^*}\\
\end{tangle}
\enspace.$ So we have
\[ (L1): \ \ \ \begin{tangle}
\object{L}\step[2]\object{L}\\
\cu*\\
\step\object{L}\\
\end{tangle}
\enspace=-\enspace
\begin{tangle}
\object{L}\step[2]\object{L}\\
\XX\\
\cu*\\
\step\object{L}\\
\end{tangle}\enspace, \ \ \  \ \ \ (L2): \ \ \
\step[1]
\begin{tangle}
\object{L}\step[2]\object{L}\step\object{L}\\
\cu*\step[1]\id\\
\step\cu*\\
\step[2]\object{L}
\end{tangle}
\enspace+\enspace
\begin{tangle}
\object{L}\step[2]\object{L}\step[2]\object{L}\\
\XX\step[2]\id\\
\id\step[2]\XX\\
\cu*\step[1]\dd\\
\step\cu*\\
\step[2]\object{L}\\
\end{tangle}
\enspace+\enspace
\begin{tangle}
\object{L}\step[2]\object{L}\step[2]\object{L}\\
\id\step[2]\XX\\
\XX\step[2]\id\\
\cu*\step[1]\dd\\
\step\cu*\\
\step[2]\object{L}\\
\end{tangle}\ \ \ , \]

\[ (CL1): \ \ \
\begin{tangle}
\step\object{L}\\
\cd*\\
\object{L}\step[2]\object{L}\\
\end{tangle}
\enspace=-\enspace
\begin{tangle}
\step[1]\object{L}\\
\cd*\\
\XX\\
\object{L}\step[2]\object{L}\\
\end{tangle}\enspace, \ \ \ \ \ (CL2):
\step[2]
\begin{tangle}
\step[2]\object{L}\\
\step\cd*\\
\cd*\step[1]\id\\
\object{L}\step[2]\object{L}\step\object{L}\\
\end{tangle}
\enspace+\enspace
\begin{tangle}
\step[2]\object{L}\\
\step\cd*\\
\cd*\step[1]\d\\
\id\step[2]\XX\\
\XX\step[2]\id\\
\object{L}\step[2]\object{L}\step[2]\object{L}\\
\end{tangle}
\enspace+\enspace
\begin{tangle}
\step[2]\object{L}\\
\step\cd*\\
\cd*\step[1]\d\\
\XX\step[2]\id\\
\id\step[2]\XX\\
\object{L}\step[2]\object{L}\step[2]\object{L}\\
\end{tangle}
\enspace =\enspace 0\enspace. \\
\]

\[ (LB) :\ \ \ \ \begin{tangle}
\object{L}\step[2]\object{L}\\
\cu*\\
\cd*\\
\object{L}\step[2]\object{L}\\
\end{tangle}
\enspace=\enspace
\begin{tangle}
\object{L}\step[3]\object{L}\\
\id\step[2]\cd*\\
\cu*\step[2]\id\\
\step\object{L}\step[3]\object{L}\\
\end{tangle}
\enspace+\enspace
\begin{tangle}
\object{L}\step[3]\object{L}\\
\id\step[2]\cd*\\
\XX\step[2]\id\\
\id\step[2]\cu*\\
\object{L}\step[3]\object{L}\\
\end{tangle}
\enspace-\enspace
\begin{tangle}
\step \object{L}\step[3]\object{L}\\
\step \XX\\
\ne1\step[1]\cd*\\
\cu*\step[2]\id\\
\step\object{L}\step[3]\object{L}\\
\end{tangle}
\enspace-\enspace
\begin{tangle}
\step\object{L}\step[3]\object{L}\\
\step\XX\\
\ne1\step[1]\cd*\\
\XX\step[2]\id\\
\id\step[2]\cu*\\
\object{L}\step[3]\object{L} \ \ \  \\
\end{tangle} \ \ \ .
\]

Obviously, if  $(L1)$ holds, then (L2) is equivalent to anyone of
the following two equations:
\[
\begin{tangle}
\object{L}\step[1]\object{L}\step[2]\object{L}\\
\id\step[1]\cu*\\
\cu*\\
\step\object{L}\\
\end{tangle}
\enspace+\enspace
\begin{tangle}
\object{L}\step[2]\object{L}\step[2]\object{L}\\
\XX\step[2]\id\\
\id\step[2]\XX\\
\d\step[1]\cu*\\
\step\cu*\\
\step[2]\object{L}\\
\end{tangle}
\enspace+\enspace
\begin{tangle}
\object{L}\step[2]\object{L}\step[2]\object{L}\\
\id\step[2]\XX\\
\XX\step[2]\id\\
\d\step[1]\cu*\\
\step\cu*\\
\step[2]\object{L}\\
\end{tangle}
\enspace =\enspace 0\enspace, \step[4]
\begin{tangle}
\object{L}\step[1]\object{L}\step[2]\object{L}\\
\id\step[1]\cu*\\
\cu*\\
\step\object{L}
\end{tangle}
\enspace=\enspace
\begin{tangle}
\object{L}\step[2]\object{L}\step\object{L}\\
\cu*\step[1]\id\\
\step\cu*\\
\step[2]\object{L}
\end{tangle}
\enspace-\enspace
\begin{tangle}
\object{L}\step[2]\object{L}\step[2]\object{L}\\
\id\step[2]\XX\\
\cu*\step[1]\dd\\
\step\cu*\\
\step[2]\object{L}
\end{tangle} \ \ \ \  .
\]

 If  $H$ is a braided Lie algebra and
   $\alpha([,] \ot \id)$ $=\alpha(\id\ot \alpha)- \alpha(\id\ot \alpha)(C_{H, H}
\otimes id)$, then $(A, \alpha )$ is called a left $H$-module. If
$(\delta\ot \id)\phi=(\id\ot \phi)\phi-(C_{H, H}\ot \id)(\id\ot
\phi)\phi$, then $(A, \phi )$ is called a left $H$-comodule. If $H$
and  $A$ are braided Lie algebra, and  $\alpha(\id\ot [$,
~$])=[,](\alpha\ot \id)+[,](\id\ot\alpha)(C _{H, A}\otimes id )$,
then $(A, [\ , \ ], \alpha )$ is called an $H$-module Lie algebra.
Assume that  $H$ is a braided Lie coalgebra and  $L$ is a braided
Lie algebra; if $L$ is a left $H$-comodule and the bracket operation
of $L$ is a morphism of  $H$-comodules, i.e.
$\phi[,]=(\id\ot[,])(\phi\ot\id)+(\id\ot[,])(C_{A, H}\ot\id)(\id\ot
\phi)$, then $L$ is called a left $H$-comodule braided Lie algebra.

Assume that $H$ is a braided Lie bialgebra,   $(A, \alpha )$ is a
left $H$-module and  $(A, \phi)$ is a left $H$-comodule. If the
following condition is satisfied:
 \[ (YD):\ \ \
\begin{tangle}
\object{H}\step[2]\object{A}\\
\lu[2]\\
\ld[2] \\
\object{H}\step[2]\object{A}
\end{tangle}
\enspace=\enspace
\begin{tangle}
\object{H}\step[4]\object{A}\\
\id\step[2]\ld[2]\\
\cu*\step[2]\id\\
\step\object{H}\step[3]\object{A}
\end{tangle}
\enspace+\enspace
\begin{tangle}
\object{H}\step[4]\object{A}\\
\id\step[2]\ld[2]\\
\XX\step[2]\id\\
\id\step[2]\lu[2]\\
\object{H}\step[4]\object{A}
\end{tangle}
\enspace+\enspace
\begin{tangle}
\step\object{H}\step[3]\object{A}\\
\cd*\step[2]\id\\
\id\step[2]\lu[2]\\
\object{H}\step[4]\object{A}
\end{tangle} \ \ \ ,
\]
then  $A$ is called a left-left- Yetter-Drinfeld module over $H$, or
(left $H$- )YD module. Denote by $\lyd H$ the category of all left
$H$-YD modules. Dually, assume that $H$ is a braided Lie bialgebra,
 $(H, \beta )$ is a right $A$-module and  $(H, \psi)$ is a right
 $A$-comodule. If the following condition is satisfied: \[ (YD):\ \ \
\begin{tangle}
\object{H}\step[2]\object{A}\\
\ru[2]\\
\rd[2] \\
\object{H}\step[2]\object{A}
\end{tangle}
\enspace=\enspace
\begin{tangle}
\object{H}\step[4]\object{A}\\

\rd[2]\step[2]\id\\
\id\step[2]\cu*\\
\object{H}\step[3]\object{A}
\end{tangle}
\enspace+\enspace
\begin{tangle}
\object{H}\step[4]\object{A}\\

\rd[2] \step[2] \id\\
  \id \step[2]  \XX\\
\ru[2]\step[2]\id\\
\object{H}\step[4]\object{A}
\end{tangle}
\enspace+\enspace
\begin{tangle}
\object{H}\step[3]\object{A}\\
\id\step[2]\cd*\\
\ru[2]\step[2]\id\\
\object{H}\step[4]\object{A}
\end{tangle} \ \ \ ,
\] then $H$ is called right-right-Yetter-Drinfeld module over $A$, or (right $A$-)YD module.
Denote by  $\ryd A$ the category of all right $A$-YD modules.

Assume that $H$ is a braided Lie bialgebra in $\cal C$, $A$ is a
left $H$-YD module and  the bracket operation and co-bracket
operation of $A$  are $H$-module morphism and $H$-comodule morphism.
If the following condition is satisfied:
\[(SLB): \ \ \
\begin{tangle}
\object{A}\step[2]\object{A}\\
\cu*\\
\cd*\\
\object{A}\step[2]\object{A}\\
\end{tangle}
\enspace=\enspace
\begin{tangle}
\object{A}\step[3]\object{A}\\
\id\step[2]\cd*\\
\cu*\step[2]\id\\
\step\object{A}\step[3]\object{A}\\
\end{tangle}
\enspace+\enspace
\begin{tangle}
\object{A}\step[3]\object{A}\\
\id\step[2]\cd*\\
\XX\step[2]\id\\
\id\step[2]\cu*\\
\object{A}\step[3]\object{A}\\
\end{tangle}
\enspace-\enspace
\begin{tangle}
\step\object{A}\step[2]\object{A}\\
\step\XX\\
\dd\step\cd*\\
\cu*\step[2]\id\\
\step\object{A}\step[3]\object{A}\\
\end{tangle}
\enspace-\enspace
\begin{tangle}
\step\object{A}\step[2]\object{A}\\
\step\XX\\
\dd\step\cd*\\
\XX\step[2]\id\\
\id\step[2]\cu*\\
\object{A}\step[3]\object{A}\\
\end{tangle}
\enspace-\enspace
\begin{tangle}
\object{A}\step[2]\object{A}\\
\ox {^S C}\\
\object{A}\step[2]\object{A}
\end{tangle} \ \ \ , \]
then $A$ is called a (left $H$-) infinitesimal braided Lie
bialgebra, where $^SC$ is called an infinitesimal braiding ( see
\cite{Ma00}) and is defined as follows:
\[
\begin{tangle}
\object{A}\step[2]\object{A}\\
\ox {^SC }\\
\object{A}\step[2]\object{A}
\end{tangle}
\enspace=\enspace
\begin{tangle}
\object{A}\step[4]\object{A}\\
\id\step[2]\Ld\\
\XX\step[2]\id\\
\Lu\step[2]\id\\
\step[2]\object{A}\step[2]\object{A}
\end{tangle}
\enspace+\enspace
\begin{tangle}
\step[2]\object{A}\step[2]\object{A}\\
\Ld\step[2]\id\\
\XX\step[2]\id\\
\id\step[2]\Lu\\
\object{A}\step[4]\object{A}
\end{tangle}
\enspace-\enspace
\begin{tangle}
\step[2]\object{A}\step[2]\object{A}\\
\Ld\step[2]\id\\
\id\step[2]\XX\\
\Lu\step[2]\id\\
\step[2]\object{A}\step[2]\object{A}
\end{tangle}
\enspace-\enspace
\begin{tangle}
\object{A}\step[4]\object{A}\\
\id\step[2]\Ld\\
\XX\step[2]\id\\
\id\step[2]\XX\\
\XX\step[2]\id\\
\id\step[2]\Lu\\
\object{A}\step[4]\object{A}
\end{tangle} \ \ \ .
\]
Similarly, we have right infinitesimal braided Lie bialgebra.

 \vspace*{.5cm}


\section {The double bicrosssum of braided Lie algebras} \label {s2}

In this section we give the conditions for the double bicrosssum of
braided Lie algebras to become a braided Lie bialgebra.

Let  $D := A \oplus H$ as objects in ${\cal C}$; define the bracket
and co-bracket  operations of $D$ as follows:
\[
[\ \ ]_D  =( \enspace
\begin{tangle}
\object{A}\step[2]\object{A}\\
\cu*\\
\step\object{A}\\
\end{tangle}
\enspace+\enspace
\begin{tangle}
\object{H}\step[2]\object{A}\\
\lu[2]\\
\step[2]\object{A}\\
\end{tangle}
\enspace-\enspace
\begin{tangle}
\object{A}\step[2]\object{H}\\
\XX\\
\lu[2]\\
\step[2]\object{A}\\
\end{tangle}
\enspace ) \oplus (\enspace
\begin{tangle}
\object{H}\step[2]\object{H}\\
\cu*\\
\step\object{H}\\
\end{tangle}
\enspace+\enspace
\begin{tangle}
\object{H}\step[2]\object{A}\\
\Ru\\
\object{H}\\
\end{tangle}
\enspace-\enspace
\begin{tangle}
\object{A}\step[2]\object{H}\\
\XX\\
\Ru\\
\object{H}\\
\end{tangle} \ \ )\ \  \ ,
\]

\[ \delta _D
= \enspace
\begin{tangle}
\step\object{A}\\
\cd*\\
\object{A}\step[2]\object{A}\\
\end{tangle}
\enspace+\enspace
\begin{tangle}
\step[2]\object{A}\\
\Ld\\
\object{H}\step[2]\object{A}\\
\end{tangle}
\enspace-\enspace
\begin{tangle}
\step[2]\object{A}\\
\Ld\\
\XX\\
\object{A}\step[2]\object{H}\\
\end{tangle}  \ \ + \ \
\begin{tangle}
\step\object{H}\\
\cd*\\
\object{H}\step[2]\object{H}\\
\end{tangle}
\enspace+\enspace
\begin{tangle}
\object{H}\\
\Rd\\
\object{H}\step[2]\object{A}\\
\end{tangle}
\enspace-\enspace
\begin{tangle}
\object{H}\\
\Rd\\
\XX\\
\object{A}\step[2]\object{H}
\end{tangle}\ \ \ .
\]
  Write  $A ^ \phi _\alpha \bowtie ^\psi _\beta H$ = $(D, [\ \ ]_D,
\delta _D)$.  If one of $\alpha, \beta, \phi$, $\psi$ is zero, we
can omit the zero morphism in $A ^ \phi _\alpha \bowtie ^\psi _\beta
H$, For example, when $\phi =0, \psi =0$ , we denote $A ^ \phi
_\alpha \bowtie ^\psi _\beta H$ by $A _\alpha \bowtie _\beta H$;
when $\alpha =0, \beta =0$, we denote $A ^ \phi _\alpha \bowtie
^\psi _\beta H$  by $A ^ \phi \bowtie ^\psi H$; Similarly, we have
 $A ^ \phi _\alpha \bowtie
 H$, $A  \bowtie ^\psi _\beta H$, $A ^ \phi
\bowtie  _\beta H$, $A  _\alpha \bowtie ^\psi H$. We call $A ^ \phi
_\alpha \bowtie ^\psi _\beta H$ the  double bicrosssum; $A _\alpha
\bowtie _\beta H$ bicrosssum; $A ^ \phi \bowtie ^\psi H$
bicrosscosum; $A ^ \phi _\alpha \bowtie
 H$ and  $A  \bowtie ^\psi _\beta H$ bisum; $A ^ \phi
\bowtie  _\beta H$ and  $A  _\alpha \bowtie ^\psi H$ bicross sum; $A
_\alpha \bowtie H$ and $A \bowtie _\beta H$ semi-direct sum, or the
smash sum; $A ^ \phi  \bowtie H$ and  $A \bowtie ^\psi H$
semi-direct co-sum. We also use the following notations: $A ^ \phi
_\alpha \bowtie ^\psi _\beta H$= $A\lrbiprod  H$; $A _\alpha \bowtie
_\beta H =A \bowtie
 H$;
 $A ^ \phi \bowtie ^\psi H$ $= A \lrcoprod H     $;
 $A
_\alpha ^\phi \bowtie
 H $ $=A\lbiprod H$    and $A  \bowtie ^\psi _\beta H$ $= A\rbiprod H$;
 $A ^ \phi
\bowtie  _\beta H     $= $A\lcprp H$ and  $A  _\alpha \bowtie ^\psi
H$=$A\lprcp H$; $A _\alpha \bowtie H$=  $A\lprod H$ and $A \bowtie
_\beta H$= $A\rprod H$; $A ^ \phi  \bowtie H$= $A\lcoprod H$ and  $A
\bowtie ^\psi H$ = $A\rcoprod H$. For conveniency, write
$$[\ \ ]_D
= \ \ \ \ \begin{tangle}
\object{A\bowtie H}\step[6]\object{A\bowtie H}\\
\hstep \step [1.5]\tu {\bowtie }\\
\step[3]\object{A\bowtie H}\\
\end{tangle} \ \ \  \mbox { and   } \ \ \ \delta _D = \ \ \ \ \  \begin{tangle}
\step[3]\object{A \lrcoprod  H}\\
\hstep \step [1.5]\td \lrcoprod\\
\object{A \lrcoprod  H}\step[5]\object{A \lrcoprod  H}

\end{tangle} \ \ \ \ \ \ \  .$$

\begin {Definition} \label {2.4} (\cite [Def.8.3.1]{Ma95}) Assume that  $A$ and  $H$ are braided Lie
algebras. If   $(A, \alpha )$ is a left $H$-module,  $(H, \beta )$
is a right $A$-module, and the following  (M1) and  (M2) hold, then
 $(A, H, \alpha, \beta )$ ( or  $(A, H )$) is called  a matched pair of braided Lie algebras:
\[ (M1): \ \ \
\begin{tangle}
\object{H}\step\object{A}\step[2]\object{A}\\
\id\step\cu*\\
\Lu\\
\step[2]\object{A}
\end{tangle}
\enspace =\enspace
\begin{tangle}
\object{H}\step[2]\object{A}\step[2]\object{A}\\
\Lu\step[2]\id\\
\step[2]\cu*\\
\step[3]\object{A}
\end{tangle}
\enspace +\enspace
\begin{tangle}
\object{H}\step[2]\object{A}\step[2]\object{A}\\
\XX\step[2]\id\\
\id\step[2]\Lu\\
\Cu*\\
\step[2]\object{A}
\end{tangle}
\enspace +\enspace
\begin{tangle}
\object{H}\step[2]\object{A}\step\object{A}\\
\Ru\step\id\\
\lu[3]\\
\step[3]\object{A}
\end{tangle}
\enspace -\enspace
\begin{tangle}
\object{H}\step[2]\object{A}\step[2]\object{A}\\
\id\step[2]\XX\\
\ru[2]\step[2]\id\\
\lu[4]\\
\step[4]\object{A}\\
\end{tangle} \ \ ; \]
\[
(M2):\ \ \
\begin{tangle}
\object{H}\step[2]\object{H}\step\object{A}\\
\cu*\step\id\\
\step\Ru\\
\step\object{H}
\end{tangle}
\enspace =\enspace
\begin{tangle}
\object{H}\step[2]\object{H}\step[2]\object{A}\\
\id\step[2]\Ru\\
\cu*\\
\step\object{H}
\end{tangle}
\enspace +\enspace
\begin{tangle}
\object{H}\step[2]\object{H}\step[2]\object{A}\\
\id\step[2]\XX\\
\Ru\step[2]\id\\
\Cu*\\
\step[2]\object{H}\\
\end{tangle}
\enspace +\enspace
\begin{tangle}
\object{H}\step\object{H}\step[2]\object{A}\\
\id\step\Lu\\
\ru[3]\\
\object{H}
\end{tangle}
\enspace -\enspace
\begin{tangle}
\object{H}\step[2]\object{H}\step[2]\object{A}\\
\XX\step[2]\id\\
\id\step[2]\Lu\\
\ru[4]\\
\object{H}
\end{tangle}\ \ \ \ .
\]
\end {Definition}

\begin {Proposition} \label {lem2.1}(\cite [Pro.8.3.2]{Ma95})
 If  $(A, H)$ is a matched pair of braided Lie algebras, then double cross sum  $A \bowtie H$ is a braided Lie algebra in
 ${\cal C}$.
\end {Proposition}

{\bf Proof.} Obviously,  (L1) hold. Now we show  (L2). By
Definition,
\[
\begin{tangle}
\object{H}\step[2]\object{A}\\
\tu \bowtie\\
\end{tangle}
\enspace =\enspace
\begin{tangle}
\object{H}\step[2]\object{A}\\
\lu[2]\\
\step[2]\object{A}
\end{tangle}
\enspace +\enspace
\begin{tangle}
\object{H}\step[2]\object{A}\\
\Ru\\
\object{H}
\end{tangle}\enspace\text{(I)} \ , \ \step[4]
\begin{tangle}
\object{A}\step[2]\object{H}\\
\tu \bowtie\\
\end{tangle}
\enspace =-\enspace
\begin{tangle}
\object{A}\step[2]\object{H}\\
\XX\\
\lu[2]\\
\step[2]\object{A}
\end{tangle}
\enspace -\enspace
\begin{tangle}
\object{A}\step[2]\object{H}\\
\XX\\
\Ru\\
\object{H}
\end{tangle}\enspace \text{(II)} \ \ \ . \]
Therefore, it follows from (I) and the conditions of matched pair
that
\begin{equation*}
\begin{tangle}
\object{H}\step\object{A}\step[2]\object{A}\\
\id\step\tu \bowtie\\
\tu \bowtie\\
\end{tangle}
\enspace = \enspace
\begin{tangle}
\object{H}\step[2]\object{A}\step[2]\object{A}\\
\lu[2]\step[2]\id\\
\step[2]\cu*\\
\step[3]\object{A}
\end{tangle}
\enspace +\enspace
\begin{tangle}
\object{H}\step[2]\object{A}\step[2]\object{A}\\
\XX\step[2]\id\\
\id\step[2]\lu[2]\\
\Cu*\\
\step[2]\object{A}
\end{tangle}
\enspace +\enspace
\begin{tangle}
\object{H}\step[2]\object{A}\step\object{A}\\
\Ru\step\id\\
\lu[3]\\
\step[3]\object{A}
\end{tangle}
\enspace -\enspace
\begin{tangle}
\object{H}\step[2]\object{A}\step[2]\object{A}\\
\id\step[2]\XX\\
\ru[2]\step[2]\id\\
\lu[4]\\
\step[4]\object{A}\\
\end{tangle}
\enspace + \enspace
\begin{tangle}
\object{H}\step[2]\object{A}\step[2]\object{A}\\
\Ru\step[2]\id\\
\ru[4]\\
\object{H}
\end{tangle}
\enspace-\enspace
\begin{tangle}
\object{H}\step[2]\object{A}\step[2]\object{A}\\
\id\step[2]\XX\\
\Ru\step\step\id\\
\ru[4]\\
\object{H}
\end{tangle} \ \ \ .
\end{equation*}
Write the six terms of the  right hand side above as (1), (2),
$\cdots$, (6), respectively.
\[(1)+(3)+(5)=
\begin{tangle}
\object{H}\step[2]\object{A}\step[2]\object{A}\\
\lu[2]\step[2]\id\\
\step[2]\cu*\\
\step[3]\object{A}
\end{tangle}
\enspace +\enspace
\begin{tangle}
\object{H}\step[2]\object{A}\step\object{A}\\
\Ru\step\id\\
\lu[3]\\
\step[3]\object{A}
\end{tangle}
\enspace +\enspace
\begin{tangle}
\object{H}\step[2]\object{A}\step[2]\object{A}\\
\Ru\step[2]\id\\
\ru[4]\\
\object{H}
\end{tangle}
\enspace =\enspace
\begin{tangle}
\object{H}\step[2]\object{A}\step\object{A}\\
\tu \bowtie\step\id\\
\step\tu \bowtie\\
\end{tangle} \ \ \ \ .
\]
Using the braided anti- symmetry of $A$, we have
\[
\begin{tangle}
\object{H}\step[2]\object{A}\step[2]\object{A}\\
\XX\step[2]\id\\
\id\step[2]\lu[2]\\
\Cu*\\
\step[2]\object{A}
\end{tangle}
\enspace -\enspace
\begin{tangle}
\object{H}\step[2]\object{A}\step[2]\object{A}\\
\id\step[2]\XX\\
\ru[2]\step[2]\id\\
\lu[4]\\
\step[4]\object{A}\\
\end{tangle}
\enspace-\enspace
\begin{tangle}
\object{H}\step[2]\object{A}\step[2]\object{A}\\
\id\step[2]\XX\\
\Ru\step\step\id\\
\ru[4]\\
\object{A}
\end{tangle}
\enspace = -\enspace
\begin{tangle}
\object{H}\step[2]\object{A}\step[2]\object{A}\\
\id\step[2]\XX\\
\lu[2]\step[2]\id\\
\step[2]\cu*\\
\step[3]\object{A}
\end{tangle}
\enspace -\enspace
\begin{tangle}
\object{H}\step[2]\object{A}\step[2]\object{A}\\
\id\step[2]\XX\\
\ru[2]\step[2]\id\\
\lu[4]\\
\step[4]\object{A}\\
\end{tangle}
\enspace-\enspace
\begin{tangle}
\object{H}\step[2]\object{A}\step[2]\object{A}\\
\id\step[2]\XX\\
\Ru\step\step\id\\
\ru[4]\\
\object{A}
\end{tangle}
\enspace =-\enspace
\begin{tangle}
\object{H}\step[2]\object{A}\step[2]\object{A}\\
\id\step[2]\XX\\
\tu \bowtie \step[1]\dd\\
\step\tu \bowtie\\
\end{tangle} \ \ \ \ .
\]
Therefore
\[
\begin{tangle}
\object{H}\step\object{A}\step[2]\object{A}\\
\id\step\tu \bowtie\\
\tu \bowtie\\
\end{tangle}
\enspace =\enspace
\begin{tangle}
\object{H}\step[2]\object{A}\step\object{A}\\
\tu \bowtie\step\id\\
\step\tu \bowtie\\
\end{tangle}
\enspace -\enspace
\begin{tangle}
\object{H}\step[2]\object{A}\step[2]\object{A}\\
\id\step[2]\XX\\
\tu \bowtie \step[1]\dd\\
\step\tu \bowtie\\
\end{tangle} \ \ \ .
\] Consequently,  (L2) holds on $H\otimes A\otimes A$. Similarly, we can show that (L2) holds for other cases.
$\Box$

In the sequel, we immediately obtain: if $(A, \alpha )$
 is an $H$-module braided Lie algebra, then semi-direct sum  $A \lprod H$ is a braided Lie
 algebra.

\begin {Definition} \label {2.5} Assume that $A$ and  $H$ are braided Lie coalgebras. If $(A,
\phi )$ is a left $H$-comodule,  $(H, \psi )$ is a right
$A$-comodule and the following (CM1) and (CM2) hold, then  $(A, H,
\phi, \psi )$  ( or  $(A, H )$)  is called a matched pair of braided
Lie coalgebras.

\[ (CM1): \ \ \ \begin{tangle}
\step[2]\object{A}\\
\Ld\\
\id\step\cd*\\
\object{H}\step\object{A}\step[2]\object{A}
\end{tangle}
\enspace =\enspace
\begin{tangle}
\step[3]\object{A}\\
\step[2]\cd*\\
\Ld\step[2]\id\\
\object{H}\step[2]\object{A}\step[2]\object{A}
\end{tangle}
\enspace +\enspace
\begin{tangle}
\step[2]\object{A}\\
\Cd*\\
\id\step[2]\Ld\\
\XX\step[2]\id\\
\object{H}\step[2]\object{A}\step[2]\object{A}
\end{tangle}
\enspace +\enspace
\begin{tangle}
\step[3]\object{A}\\
\ld[3]\\
\Rd\step\id\\
\object{H}\step[2]\object{A}\step\object{A}
\end{tangle}
\enspace -\enspace
\begin{tangle}
\step[4]\object{A}\\
\ld[4]\\
\rd[2]\step[2]\id\\
\id\step[2]\XX\\
\object{H}\step[2]\object{A}\step[2]\object{A}
\end{tangle} \ \ ; \]

\[ (CM 2): \ \ \ \begin{tangle}
\step\object{H}\\
\step\Rd\\
\cd*\step\id\\
\object{H}\step[2]\object{H}\step\object{A}
\end{tangle}
\enspace =\enspace
\begin{tangle}
\step\object{H}\\
\cd*\\
\id\step[2]\Rd\\
\object{H}\step[2]\object{H}\step[2]\object{A}
\end{tangle}
\enspace +\enspace
\begin{tangle}
\step[2]\object{H}\\
\Cd*\\
\Rd\step[2]\id\\
\id\step[2]\XX\\
\object{H}\step[2]\object{H}\step[2]\object{A}
\end{tangle}
\enspace +\enspace
\begin{tangle}
\object{H}\\
\rd[3]\\
\id\step\Ld\\
\object{H}\step\object{H}\step[2]\object{A}
\end{tangle}
\enspace -\enspace
\begin{tangle}
\object{H}\\
\rd[4]\\
\id\step[2]\Ld\\
\XX\step[2]\id\\
\object{H}\step[2]\object{H}\step[2]\object{A}
\end{tangle}\ \ \ . \]

\end  {Definition}

\begin {Proposition}\label {lem2.2}
 If $(A, H)$ is a matched pair of braided Lie coalgebras, then  double cross co-sum $A\lrcoprod H$ is a braided Lie
 coalgebra.
\end {Proposition}

{\bf Proof.} We only check the braided  co-Jacobi identity. By
definition,
\[
\begin{tangle}
\step[2]\object{A}\\
\step\td {\lrcoprod}\\
\td {\lrcoprod}\step\id\\
\end{tangle}
\enspace =  \enspace
\begin{tangle}
\step[2]\object{A}\\
\step\cd*\\
\cd*\step\id\\
\object{A}\step[2]\object{A}\step\object{A}\\
\end{tangle}
\enspace (1)+ \enspace
\begin{tangle}
\step[2]\object{A}\\
\step\cd*\\
\ld\step[2]\id\\
\object{H}\step[2]\object{A}\step\object{A}\\
\end{tangle}
\enspace(2) - \enspace
\begin{tangle}
\step[2]\object{A}\\
\step\cd*\\
\ld\step[2]\id\\
\X\step[2]\id\\
\object{A}\step[2]\object{H}\step\object{A}\\
\end{tangle}
\enspace(3) +\enspace
\begin{tangle}
\step[3]\object{A}\\
\step\Ld\\
\cd*\step\id\\
\object{H}\step[2]\object{H}\step\object{A}\\
\end{tangle}
\enspace (4)+\enspace
\begin{tangle}
\step[2]\object{A}\\
\ld[2]\\
\rd\step\id\\
\object{H}\step\object{A}\step\object{A}\\
\end{tangle}\enspace(5)
\]
\[
\enspace -\enspace
\begin{tangle}
\step[2]\object{A}\\
\ld[2]\\
\rd\step\id\\
\X\step\id\\
\object{A}\step\object{H}\step\object{A}\\
\end{tangle}
\enspace(6)  -\enspace
\begin{tangle}
\step[3]\object{A}\\
\step\Ld\\
\step\XX\\
\cd*\step\id\\
\object{A}\step[2]\object{A}\step\object{H}\\
\end{tangle}
\enspace (7) -\enspace
\begin{tangle}
\step[3]\object{A}\\
\step\Ld\\
\step\XX\\
\ld\step[2]\id\\
\object{H}\step\object{A}\step[2]\object{H}\\
\end{tangle}
\enspace (8)+\enspace
\begin{tangle}
\step[3]\object{A}\\
\step\Ld\\
\step\XX\\
\ld\step[2]\id\\
\X\step[2]\id\\
\object{A}\step\object{H}\step[2]\object{H}\ \ \ . \\
\end{tangle}\enspace (9) \ \ \ .
\]
So
\[
\begin{tangle}
\step[2]\object{A}\\
\step\td {\lrcoprod}\\
\td {\lrcoprod}\step\id\\
\id\step[2]\X\\
\XX\step\id\\
\end{tangle}
\enspace =  \enspace
\begin{tangle}
\step[2]\object{A}\\
\step\cd*\\
\cd*\step\id\\
\id\step[2]\X\\
\XX\step\id\\
\object{A}\step[2]\object{A}\step\object{A}\\
\end{tangle}
\enspace(10) + \enspace
\begin{tangle}
\step[2]\object{A}\\
\step\cd*\\
\ld\step[2]\id\\
\id\step\XX\\
\X\step[2]\id\\
\object{A}\step\object{H}\step[2]\object{A}\\
\end{tangle}
\enspace (11)-\enspace
\begin{tangle}
\step[2]\object{A}\\
\step\cd*\\
\ld\step[2]\id\\
\X\step[2]\id\\
\id\step\XX\\
\X\step[2]\id\\
\object{A}\step\object{A}\step[2]\object{H}\\
\end{tangle}
\enspace (12)+\enspace
\begin{tangle}
\step[3]\object{A}\\
\step\Ld\\
\cd*\step\id\\
\id\step[2]\X\\
\XX\step\id\\
\object{A}\step[2]\object{H}\step\object{H}\\
\end{tangle}
\enspace(13) +\enspace
\begin{tangle}
\step[2]\object{A}\\
\ld[2]\\
\rd\step\id\\
\id\step\X\\
\X\step\id\\
\object{A}\step\object{H}\step\object{A}\\
\end{tangle}\enspace(14)
\]
\[
\enspace -\enspace
\begin{tangle}
\step[2]\object{A}\\
\ld[2]\\
\rd\step\id\\
\X\step\id\\
\id\step\X\\
\X\step\id\\
\object{A}\step\object{A}\step\object{H}\\
\end{tangle}
\enspace(15) -\enspace
\begin{tangle}
\step[2]\object{A}\\
\ld[2]\\
\id\step\cd*\\
\object{H}\step\object{A}\step[2]\object{A}\\
\end{tangle}
\enspace (16)-\enspace
\begin{tangle}
\step[2]\object{A}\\
\ld[2]\\
\id\step\ld\\
\object{H}\step\object{H}\step\object{A}\\
\end{tangle}
\enspace (17)+\enspace
\begin{tangle}
\step[2]\object{A}\\
\ld[2]\\
\id\step\ld\\
\id\step\X\\
\object{H}\step\object{A}\step\object{H}\\
\end{tangle}\enspace (18) \ \ \ \mbox { and }
\]
\[
\begin{tangle}
\step[2]\object{A}\\
\step\td {\lrcoprod}\\
\td {\lrcoprod}\step\id\\
\XX\step\id\\
\id\step[2]\X\\
\end{tangle}
\enspace =  \enspace
\begin{tangle}
\step[2]\object{A}\\
\step\cd*\\
\cd*\step\id\\
\XX\step\id\\
\id\step[2]\X\\
\object{A}\step[2]\object{A}\step\object{A}\\
\end{tangle}
\enspace (19) + \enspace
\begin{tangle}
\step[2]\object{A}\\
\step\cd*\\
\ld\step[2]\id\\
\X\step[2]\id\\
\id\step\XX\\
\object{A}\step\object{A}\step[2]\object{H}\\
\end{tangle}
\enspace (20) -\enspace
\begin{tangle}
\step[2]\object{A}\\
\step\cd*\\
\ld\step[2]\id\\
\id\step\XX\\
\object{H}\step\object{A}\step[2]\object{A}\\
\end{tangle}
\enspace (21)+\enspace
\begin{tangle}
\step[3]\object{A}\\
\step\Ld\\
\cd*\step\id\\
\XX\step\id\\
\id\step[2]\X\\
\object{H}\step[2]\object{A}\step\object{H}\\
\end{tangle}
\enspace (22)+\enspace
\begin{tangle}
\step[2]\object{A}\\
\ld[2]\\
\rd\step\id\\
\X\step\id\\
\id\step\X\\
\object{A}\step\object{A}\step\object{H}\\
\end{tangle}\enspace (23)
\]
\[ \enspace -\enspace
\begin{tangle}
\step[2]\object{A}\\
\ld[2]\\
\rd\step\id\\
\id\step\X\\
\object{H}\step\object{A}\step\object{A}\\
\end{tangle}
\enspace (24)-\enspace
\begin{tangle}
\step[2]\object{A}\\
\ld[2]\\
\id\step\cd*\\
\id\step\XX\\
\X\step[2]\id\\
\object{A}\step\object{H}\step[2]\object{A}\\
\end{tangle}
\enspace (25)-\enspace
\begin{tangle}
\step[2]\object{A}\\
\ld[2]\\
\id\step\ld\\
\id\step\X\\
\X\step\id\\
\object{A}\step\object{H}\step\object{H}\\
\end{tangle}
\enspace (26)+\enspace
\begin{tangle}
\step[2]\object{A}\\
\ld[2]\\
\id\step\ld\\
\X\step\id\\
\object{H}\step\object{H}\step\object{A}\\
\end{tangle}\enspace (27) \ \ \ .
\]
It follows from the braided Jacobi identity of $A$ that
$(1)+(10)+(19)=0$; Since $A$ is a left $H$-comodule,
$(13)+(9)-(26)=0$, $(4)-(17)+(27)=0$, $(22)-(8)+(18)=0$; By the
conditions of matched pair of braided Lie coalgebras,
$-(16)+(2)-(21)+(5)-(24)=0$, $-(25)+(11)-(3)+(14)-(6)=0$,
$-(7)+(20)-(12)+(23)-(15)=0$. Therefore, (CL2) holds on $A$.
Similarly, (CL2) holds on $H$.
 $\Box$

In the sequel, we immediately obtain: if  $(A, \phi )$
 is an $H$-comodule braided Lie coalgebra, then the semi-direct sum  $A \lcoprod
H$ is a braided Lie coalgebra.

\begin {Definition} \label {2.5''}
 Assume that  $A$
and  $H$ are braided Lie algebras and braided Lie coalgebras, and
(SLB) holds on $A \otimes A$ and  $H \otimes H $. If the following
(B1)- (B5) hold, then  $(A, H, \alpha, \beta, \phi, \psi )$ ( or
$(A, H )$) is called a double matched pair:
\[ (B1): \ \ \
\begin{tangle}
\object{H}\step[2]\object{A}\\
\Lu\\
\step\cd*\\
\step\object{A}\step[2]\object{A}
\end{tangle}
\enspace =\enspace
\begin{tangle}
\object{H}\step[3]\object{A}\\
\id\step[2]\cd*\\
\Lu\step[2]\id\\
\step[2]\object{A}\step[2]\object{A}\\
\end{tangle}
\enspace +\enspace
\begin{tangle}
\object{H}\step[3]\object{A}\\
\id\step[2]\cd*\\
\XX\step[2]\id\\
\id\step[2]\Lu\\
\object{A}\step[4]\object{A}\\
\end{tangle}
\enspace +\enspace
\begin{tangle}
\object{H}\step[4]\object{A}\\
\Rd\step[2]\id\\
\id\step[2]\XX\\
\Lu\step[2]\id\\
\step[2]\object{A}\step[2]\object{A}\\
\end{tangle}
\enspace -\enspace
\begin{tangle}
\object{H}\step[4]\object{A}\\
\Rd\step[2]\id\\
\XX\step[2]\id\\
\id\step[2]\Lu\\
\object{A}\step[4]\object{A}\\
\end{tangle}
\ \ ;
\]

\[ (B2): \ \ \
\begin{tangle}
\step\object{H}\step[2]\object{A}\\
\step\Ru\\
\cd*\\
\object{H}\step[2]\object{H}
\end{tangle}
\enspace =\enspace
\begin{tangle}
\step\object{H}\step[3]\object{A}\\
\cd*\step[2]\id\\
\id\step[2]\Ru\\
\object{H}\step[2]\object{H}\\
\end{tangle}
\enspace +\enspace
\begin{tangle}
\step\object{H}\step[3]\object{A}\\
\cd*\step[2]\id\\
\id\step[2]\XX\\
\Ru\step[2]\id\\
\object{H}\step[4]\object{H}\\
\end{tangle}
\enspace +\enspace
\begin{tangle}
\object{H}\step[4]\object{A}\\
\id\step[2]\Ld\\
\XX\step[2]\id\\
\id\step[2]\Ru\\
\object{H}\step[2]\object{H}\\
\end{tangle}
\enspace -\enspace
\begin{tangle}
\object{H}\step[4]\object{A}\\
\id\step[2]\Ld\\
\id\step[2]\XX\\
\Ru\step[2]\id\\
\object{H}\step[4]\object{H}\\
\end{tangle}
\ \ ; \]

\[ (B3): \ \ \
\begin{tangle}
\step\object{A}\step[2]\object{A}\\
\step\cu*\\
\Ld\\
\object{H}\step[2]\object{A}\\
\end{tangle}
\enspace =\enspace
\begin{tangle}
\step[2]\object{A}\step[2]\object{A}\\
\Ld\step[2]\id\\
\id\step[2]\cu*\\
\object{H}\step[3]\object{A}\\
\end{tangle}
\enspace +\enspace
\begin{tangle}
\object{A}\step[4]\object{A}\\
\id\step[2]\Ld\\
\XX\step[2]\id\\
\id\step[2]\cu*\\
\object{H}\step[3]\object{A}\\
\end{tangle}
\enspace +\enspace
\begin{tangle}
\step[2]\object{A}\step[2]\object{A}\\
\Ld\step[2]\id\\
\id\step[2]\XX\\
\Ru\step[2]\id\\
\object{H}\step[4]\object{A}\\
\end{tangle}
\enspace -\enspace
\begin{tangle}
\object{A}\step[4]\object{A}\\
\id\step[2]\Ld\\
\XX\step[2]\id\\
\Ru\step[2]\id\\
\object{H}\step[4]\object{A}\\
\end{tangle}
\ \ ; \]

\[ (B4): \ \ \
\begin{tangle}
\object{H}\step[2]\object{H}\\
\cu*\\
\step\Rd\\
\step\object{H}\step[2]\object{A}\\
\end{tangle}
\enspace =\enspace
\begin{tangle}
\object{H}\step[2]\object{H}\\
\id\step[2]\Rd\\
\cu*\step[2]\id\\
\step\object{H}\step[3]\object{A}\\
\end{tangle}
\enspace +\enspace
\begin{tangle}
\object{H}\step[4]\object{H}\\
\Rd\step[2]\id\\
\id\step[2]\XX\\
\cu*\step[2]\id\\
\step\object{H}\step[3]\object{A}\\
\end{tangle}
\enspace +\enspace
\begin{tangle}
\object{H}\step[2]\object{H}\\
\id\step[2]\Rd\\
\XX\step[2]\id\\
\id\step[2]\Lu\\
\object{H}\step[4]\object{A}\\
\end{tangle}
\enspace -\enspace
\begin{tangle}
\object{H}\step[4]\object{H}\\
\Rd\step[2]\id\\
\id\step[2]\XX\\
\id\step[2]\Lu\\
\object{H}\step[4]\object{A}\\
\end{tangle}
\ \ ; \]

 (B5):  \[ \begin{tangle}
\object{H}\step[2]\object{A}\\
\lu[2]\\
\ld[2] \\
\object{H}\step[2]\object{A}
\end{tangle}\ + \
\begin{tangle}
\object{H}\step[2]\object{A}\\
\ru[2]\\
\rd[2] \\
\object{H}\step[2]\object{A}
\end{tangle}
\enspace=\enspace
\begin{tangle}
\object{H}\step[4]\object{A}\\
\id\step[2]\ld[2]\\
\cu*\step[2]\id\\
\step\object{H}\step[3]\object{A}
\end{tangle}
\enspace+\enspace
\begin{tangle}
\object{H}\step[4]\object{A}\\
\id\step[2]\ld[2]\\
\XX\step[2]\id\\
\id\step[2]\lu[2]\\
\object{H}\step[4]\object{A}
\end{tangle}
\enspace+\enspace
\begin{tangle}
\step\object{H}\step[3]\object{A}\\
\cd*\step[2]\id\\
\id\step[2]\lu[2]\\
\object{H}\step[4]\object{A}
\end{tangle} \enspace+\enspace
\begin{tangle}
\object{H}\step[4]\object{A}\\

\rd[2]\step[2]\id\\
\id\step[2]\cu*\\
\object{H}\step[3]\object{A}
\end{tangle}
\enspace+\enspace
\begin{tangle}
\object{H}\step[4]\object{A}\\

\rd[2] \step[2] \id\\
  \id \step[2]  \XX\\
\ru[2]\step[2]\id\\
\object{H}\step[4]\object{A}
\end{tangle}
\enspace+\enspace
\begin{tangle}
\object{H}\step[3]\object{A}\\
\id\step[2]\cd*\\
\ru[2]\step[2]\id\\
\object{H}\step[4]\object{A}
\end{tangle} \ \ \ .
\]
\end {Definition}
\begin {Theorem} \label{thmain}
 If  $(A,H)$ is  matched pairs of braided Lie algebras and braided Lie
 coalgebras in $\cal
C$, and a double matched pair, then  $A\lrbiprod H$ is a braided Lie
bialgebra in $\cal C$. \end  {Theorem}

{\bf Proof. } By Proposition  \ref {lem2.1} and Proposition \ref
{lem2.2}, $A\lrbiprod H$ is a braided Lie algebra and braided Lie
coalgebra. Therefore, we have  to check
 (LB). By definition,
\begin{equation*}
\begin{tangle}
\object{H}\step[2]\object{A}\\
\tu {\bowtie}\\
\td {\lrcoprod}\\
\end{tangle}
\enspace = \enspace
\begin{tangle}
\object{H}\step[2]\object{A}\\
\Lu\\
\step\cd*\\
\step\object{A}\step[2]\object{A}\\
\end{tangle}
\enspace(a)+\enspace
\begin{tangle}
\object{H}\step[2]\object{A}\\
\Lu\\
\Ld\\
\object{H}\step[2]\object{A}\\
\end{tangle}
\enspace(b)-\enspace
\begin{tangle}
\object{H}\step[2]\object{A}\\
\Lu\\
\Ld\\
\XX\\
\object{A}\step[2]\object{H}\\
\end{tangle}
\enspace(c) +\enspace
\begin{tangle}
\step\object{H}\step[2]\object{A}\\
\step\Ru\\
\cd*\\
\object{H}\step[2]\object{H}\\
\end{tangle}
\enspace(d)+\enspace
\begin{tangle}
\object{H}\step[2]\object{A}\\
\Ru\\
\Rd\\
\object{H}\step[2]\object{A}\\
\end{tangle}
\enspace(e)-\enspace
\begin{tangle}
\object{H}\step[2]\object{A}\\
\Ru\\
\Rd\\
\XX\\
\object{A}\step[2]\object{H}\\
\end{tangle}\enspace(f) \ \ \ .
\end{equation*}
Using (B1), (B2) and (B5) on the right hand side above, we have
\[
\begin{tangle}
\object{H}\step[2]\object{A}\\
\tu {\bowtie}\\
\td {\lrcoprod}\\
\end{tangle}
\enspace = \enspace
\begin{tangle}
\object{H}\step[2]\object{A}\\
\id\step\cd*\\
\lu\step[2]\id\\
\step\object{A}\step[2]\object{A}
\end{tangle}
\enspace(1) +\enspace
\begin{tangle}
\object{H}\step[2]\object{A}\\
\id\step\cd*\\
\X\step[2]\id\\
\id\step\Lu\\
\object{A}\step[3]\object{A}
\end{tangle}
\enspace(2) +\enspace
\begin{tangle}
\object{H}\step[3]\object{A}\\
\Rd\step\id\\
\id\step[2]\X\\
\Lu\step\id\\
\step[2]\object{A}\step\object{A}\\
\end{tangle}\enspace(3)
-\enspace
\begin{tangle}
\object{H}\step[3]\object{A}\\
\Rd\step\id\\
\XX\step\id\\
\id\step[2]\lu\\
\object{A}\step[3]\object{A}\\
\end{tangle}
\enspace(4)\step[4]
\]
\[+\enspace
\begin{tangle}
\object{H}\step[3]\object{A}\\
\id\step[2]\ld\\
\cu*\step\id\\
\step\object{H}\step[2]\object{A}
\end{tangle}
\enspace(5)+\enspace
\begin{tangle}
\object{H}\step[3]\object{A}\\
\id\step[2]\ld\\
\XX\step\id\\
\id\step[2]\lu\\
\object{H}\step[3]\object{A}
\end{tangle}
\enspace(6)+\enspace
\begin{tangle}
\step\object{H}\step[2]\object{A}\\
\cd*\step\id\\
\id\step[2]\lu\\
\object{H}\step[3]\object{A}
\end{tangle}
\enspace(7)-\enspace
\begin{tangle}
\object{H}\step[3]\object{A}\\
\id\step[2]\ld\\
\cu*\step\id\\
\step\XX\\
\step\object{A}\step[2]\object{H}
\end{tangle}
\enspace(8)
\]
\[-\enspace
\begin{tangle}
\object{H}\step[3]\object{A}\\
\id\step[2]\ld\\
\XX\step\id\\
\d\step\lu\\
\step\XX\\
\step \object{A}\step[2]\object{H}
\end{tangle}
\enspace(9)-\enspace
\begin{tangle}
\step\object{H}\step[2]\object{A}\\
\cd*\step\id\\
\d\step\lu\\
\step\XX\\
\step\object{A}\step[2]\object{H}
\end{tangle}
\enspace(10) +\enspace
\begin{tangle}
\step\object{H}\step[2]\object{A}\\
\cd*\step\id\\
\id\step[2]\ru\\
\object{H}\step[2]\object{H}
\end{tangle}
\enspace(11) +\enspace
\begin{tangle}
\step\object{H}\step[2]\object{A}\\
\cd*\step\id\\
\id\step[2]\X\\
\Ru\step\id\\
\object{H}\step[3]\object{H}
\end{tangle}
\enspace(12)
\]
\[ +\enspace
\begin{tangle}
\object{H}\step[3]\object{A}\\
\id\step[2]\ld\\
\XX\step\id\\
\id\step[2]\ru\\
\object{H}\step[2]\object{H}\\
\end{tangle}\enspace(13)
 -\enspace
\begin{tangle}
\object{H}\step[3]\object{A}\\
\id\step\Ld\\
\id\step\XX\\
\ru\step[2]\id\\
\object{H}\step[3]\object{H}\\
\end{tangle}
\enspace(14) +\enspace
\begin{tangle}
\object{H}\step[3]\object{A}\\
\rd\step[2]\id\\
\id\step\cu*\\
\object{H}\step[2]\object{A}
\end{tangle}
\enspace(15) +\enspace
\begin{tangle}
\object{H}\step[3]\object{A}\\
\rd\step[2]\id\\
\id\step\XX\\
\ru\step[2]\id\\
\object{H}\step[3]\object{A}
\end{tangle}
\enspace(16)
\]
\[
+\enspace
\begin{tangle}
\object{H}\step[2]\object{A}\\
\id\step\cd*\\
\ru\step[2]\id\\
\object{H}\step[3]\object{A}
\end{tangle}
\enspace(17) -\enspace
\begin{tangle}
\object{H}\step[3]\object{A}\\
\rd\step[2]\id\\
\id\step\cu*\\
\XX\\
\object{A}\step[2]\object{H}
\end{tangle}
\enspace(18)- \enspace
\begin{tangle}
\object{H}\step[3]\object{A}\\
\rd\step[2]\id\\
\id\step\XX\\

\ru\step[1]\ne1\\

\XX\\
\object{A}\step[2]\object{H}
\end{tangle}
\enspace(19) -\enspace
\begin{tangle}
\object{H}\step[2]\object{A}\\
\id\step\cd*\\
\ru\step\dd\\
\XX\\
\object{A}\step[2]\object{H}
\end{tangle}
\enspace(20) \ \ \ .
\]
 We have to show the above is equal to
\[
\begin{tangle}
\object{H}\step[3]\object{A}\\
\id\step[2]\td \lrcoprod\\
\tu \bowtie\step[2]\id\\
\end{tangle}
\enspace(g)+\enspace
\begin{tangle}
\object{H}\step[3]\object{A}\\
\id\step[2]\td \lrcoprod\\
\XX\step[2]\id\\
\id\step[2]\tu \bowtie\\
\end{tangle}
\enspace(h)-\enspace
\begin{tangle}
\step\object{H}\step[2]\object{A}\\
\step\XX\\
\dd\step[1]\td \lrcoprod\\
\tu \bowtie\step[2]\id\\
\end{tangle}
\enspace(i)-\enspace
\begin{tangle}
\step\object{H}\step[2]\object{A}\\
\step\XX\\
\dd\step[1]\td \lrcoprod\\
\XX\step[2]\id\\
\id\step[2]\tu \bowtie\\
\end{tangle}\enspace(j) \ \ \ .
\]
Indeed,
\[
(g)=\enspace
\begin{tangle}
\object{H}\step[2]\object{A}\\
\id\step\cd*\\
\lu\step[2]\id\\
\step\object{A}\step[2]\object{A}
\end{tangle}
\enspace(1) +\enspace
\begin{tangle}
\object{H}\step[2]\object{A}\\
\id\step\cd*\\
\ru\step[2]\id\\
\object{H}\step[3]\object{A}
\end{tangle}
\enspace(17) +\enspace
\begin{tangle}
\object{H}\step[3]\object{A}\\
\id\step[2]\ld\\
\cu*\step\id\\
\step\object{H}\step[2]\object{A}
\end{tangle}
\enspace(5)  -\enspace
\begin{tangle}
\object{H}\step[3]\object{A}\\
\id\step\Ld\\
\id\step\XX\\
\lu\step[2]\id\\
\step\object{A}\step[2]\object{H}
\end{tangle}
\enspace(9') -\enspace
\begin{tangle}
\object{H}\step[3]\object{A}\\
\id\step\Ld\\
\id\step\XX\\
\ru\step[2]\id\\
\object{H}\step[3]\object{H}\\
\end{tangle}
\enspace(14) \ \ \ ,
\]
\[
(h)=\enspace
\begin{tangle}
\object{H}\step[2]\object{A}\\
\id\step\cd*\\
\X\step[2]\id\\
\id\step\Lu\\
\object{A}\step[3]\object{A}
\end{tangle}
\enspace(2) +\enspace
\begin{tangle}
\object{H}\step[2]\object{A}\\
\id\step\cd*\\
\X\step[2]\id\\
\id\step\Ru\\
\object{A}\step[2]\object{H}
\end{tangle}
\enspace(20')+\enspace
\begin{tangle}
\object{H}\step[3]\object{A}\\
\id\step[2]\ld\\
\XX\step\id\\
\id\step[2]\lu\\
\object{H}\step[3]\object{A}
\end{tangle}
\enspace(6)+\enspace
\begin{tangle}
\object{H}\step[3]\object{A}\\
\id\step[2]\ld\\
\XX\step\id\\
\id\step[2]\ru\\
\object{H}\step[2]\object{H}\\
\end{tangle}\enspace(13)-\enspace
\begin{tangle}
\object{H}\step[3]\object{A}\\
\id\step[2]\ld\\
\cu*\step\id\\
\step\XX\\
\step\object{A}\step[2]\object{H}
\end{tangle}
\enspace(8) \ \ \ ,
\]

\[
-(i)=\enspace
\begin{tangle}
\object{H}\step[2]\object{A}\\
\XX\\
\id\step\cd*\\
\X\step[2]\id\\
\lu\step[2]\id\\
\step\object{A}\step[2]\object{H}
\end{tangle}
\enspace(10') +
\begin{tangle}
\object{H}\step[2]\object{A}\\
\XX\\
\id\step\cd*\\
\X\step[2]\id\\
\ru\step[2]\id\\
\object{H}\step[3]\object{H}
\end{tangle}
\enspace(12') +\enspace
\begin{tangle}
\object{H}\step[3]\object{A}\\
\XX\\
\id\step[2]\rd\\
\XX\step\id\\
\Lu\step\id\\
\step[2]\object{A}\step\object{H}
\end{tangle}
\enspace(3') +\enspace
\begin{tangle}
\object{H}\step[2]\object{A}\\
\XX\\
\id\step[2]\rd\\
\XX\step\id\\
\Ru\step\id\\
\object{A}\step[2]\object{H}
\end{tangle}
\enspace(16')+
\begin{tangle}
\object{H}\step[2]\object{A}\\
\XX\\
\id\step[2]\rd\\
\id\step[2]\X\\
\cu*\step\id\\
\step\object{A}\step[2]\object{H}\\
\end{tangle}
\enspace(18') \ \ \ ,
\]
\[
-(j)=\enspace
\begin{tangle}
\object{H}\step[2]\object{A}\\
\XX\\
\id\step\cd*\\
\X\step[2]\id\\
\id\step\XX\\
\id\step\Lu\\
\object{H}\step[3]\object{A}
\end{tangle}
\enspace(7') +
\begin{tangle}
\object{H}\step[2]\object{A}\\
\XX\\
\id\step\cd*\\
\X\step[2]\id\\
\id\step\XX\\
\id\step\Ru\\
\object{H}\step[2]\object{H}
\end{tangle}
\enspace(11')-\enspace
\begin{tangle}
\object{H}\step[2]\object{A}\\
\XX\\
\id\step[2]\Rd\\
\XX\step[2]\id\\
\id\step[2]\cu*\\
\object{H}\step[3]\object{A}
\end{tangle}
\enspace(15')-\enspace
\begin{tangle}
\object{H}\step[2]\object{A}\\
\XX\\
\id\step[2]\rd\\
\id\step[2]\X\\
\XX\step\id\\
\id\step[2]\X\\
\id\step[2]\lu\\
\object{H}\step[3]\object{H}\\
\end{tangle}\enspace(4')-
\begin{tangle}
\object{H}\step[2]\object{A}\\
\XX\\
\id\step[2]\rd\\
\id\step[2]\X\\
\XX\step\id\\
\id\step[2]\X\\
\id\step[2]\ru\\
\object{H}\step[2]\object{H}\\
\end{tangle}\enspace(19') \ \ \ .
\]
Thus  (LB) holds on  $H \otimes A$. Similarly, it follows from (B1),
(B2) and  (B5) that  (LB) holds on  $A\otimes H$.

Next we investigate the case  of (LB) on   $A \otimes A$:
\[
\begin{tangle}
\object{A}\step[2]\object{A}\\
\tu \lprod\\
\td \lrcoprod\\
\end{tangle}
\enspace = \enspace
\begin{tangle}
\object{A}\step[2]\object{A}\\
\cu*\\
\cd*\\
\object{A}\step[2]\object{A}\\
\end{tangle}
\enspace(k)+\enspace
\begin{tangle}
\step\object{A}\step[2]\object{A}\\
\step\cu*\\
\Ld\\
\object{H}\step[2]\object{A}\\
\end{tangle}
\enspace(l)-\enspace
\begin{tangle}
\step\object{A}\step[2]\object{A}\\
\step\cu*\\
\Ld\\
\XX\\
\object{A}\step[2]\object{H}\\
\end{tangle}\enspace(m) \ \ \ ,
\]
where
\[
(k)  \stackrel {(SLB)}{=} \enspace
\begin{tangle}
\object{A}\step[3]\object{A}\\
\id\step[2]\cd*\\
\cu*\step[2]\id\\
\step\object{A}\step[3]\object{A}\\
\end{tangle}
\enspace(1)-\enspace
\begin{tangle}
\step\object{A}\step[2]\object{A}\\
\step\XX\\
\dd\step[1]\cd*\\
\cu*\step[2]\id\\
\step\object{A}\step[3]\object{A}\\
\end{tangle}
\enspace(2)+\enspace
\begin{tangle}
\object{A}\step[3]\object{A}\\
\id\step[2]\cd*\\
\XX\step[2]\id\\
\id\step[2]\cu*\\
\object{A}\step[3]\object{A}\\
\end{tangle}
\enspace(3)-\enspace
\begin{tangle}
\step\object{A}\step[2]\object{A}\\
\step\XX\\
\dd\step[1]\cd*\\
\XX\step[2]\id\\
\id\step[2]\cu*\\
\object{A}\step[3]\object{A}\\
\end{tangle}\enspace(4)
\]
\[
\enspace-\enspace
\begin{tangle}
\object{A}\step[4]\object{A}\\
\id\step[2]\Ld\\
\XX\step[2]\id\\
\Lu\step[2]\id\\
\step[2]\object{A}\step[2]\object{A}
\end{tangle}
\enspace(5)-\enspace
\begin{tangle}
\step[2]\object{A}\step[2]\object{A}\\
\Ld\step[2]\id\\
\XX\step[2]\id\\
\id\step[2]\Lu\\
\object{A}\step[4]\object{A}
\end{tangle}
\enspace(6)+\enspace
\begin{tangle}
\step[2]\object{A}\step[2]\object{A}\\
\Ld\step[2]\id\\
\id\step[2]\XX\\
\Lu\step[2]\id\\
\step[2]\object{A}\step[2]\object{A}
\end{tangle}
\enspace(7)+\enspace
\begin{tangle}
\object{A}\step[4]\object{A}\\
\id\step[2]\Ld\\
\XX\step[2]\id\\
\id\step[2]\XX\\
\XX\step[2]\id\\
\id\step[2]\Lu\\
\object{A}\step[4]\object{A}
\end{tangle}
\enspace(8) \ \ \ \ ,
\]
\[
(l)\stackrel {(B3)}{=}\enspace
\begin{tangle}
\step[2]\object{A}\step[2]\object{A}\\
\Ld\step[2]\id\\
\id\step[2]\cu*\\
\object{H}\step[3]\object{A}\\
\end{tangle}
\enspace (9)+\enspace
\begin{tangle}
\object{A}\step[4]\object{A}\\
\id\step[2]\Ld\\
\XX\step[2]\id\\
\id\step[2]\cu*\\
\object{H}\step[3]\object{A}\\
\end{tangle}
\enspace(10) +\enspace
\begin{tangle}
\step[2]\object{A}\step[2]\object{A}\\
\Ld\step[2]\id\\
\id\step[2]\XX\\
\Ru\step[2]\id\\
\object{H}\step[4]\object{A}\\
\end{tangle}
\enspace(11) -\enspace
\begin{tangle}
\object{A}\step[4]\object{A}\\
\id\step[2]\Ld\\
\XX\step[2]\id\\
\Ru\step[2]\id\\
\object{H}\step[4]\object{A}\\
\end{tangle}\enspace(12) \ \ \ \ ,
\]
\[
-(m)=-\enspace
\begin{tangle}
\step\object{A}\step[2]\object{A}\\
\ld\step[2]\id\\
\id\step\cu*\\
\XX\\
\object{A}\step[2]\object{H}\\
\end{tangle}
\enspace(13)-\enspace
\begin{tangle}
\object{A}\step[4]\object{A}\\
\id\step[2]\Ld\\
\id\step[2]\XX\\
\cu*\step[2]\id\\
\step\object{A}\step[3]\object{H}\\
\end{tangle}
\enspace(14)-\enspace
\begin{tangle}
\step[2]\object{A}\step[2]\object{A}\\
\Ld\step[2]\id\\
\XX\step[2]\id\\
\id\step[2]\Ru\\
\object{A}\step[2]\object{H}\\
\end{tangle}
\enspace(15)+\enspace
\begin{tangle}
\object{A}\step[3]\object{A}\\
\id\step[2]\ld\\
\XX\step[1]\id\\
\Ru\dd\\
\XX\\
\object{A}\step[2]\object{H}\\
\end{tangle}\enspace(16) \ \ \ \ .
\]

\mbox {the right hand side  of } (LB) =
\[ \begin{tangle}
\object{A}\step[3]\object{A}\\
\id\step[2]\td \lrcoprod\\
\tu \bowtie\step[2]\id\\
\end{tangle}
\enspace(n)+\enspace
\begin{tangle}
\object{A}\step[3]\object{A}\\
\id\step[2]\td \lrcoprod\\
\XX\step[2]\id\\
\id\step[2]\tu \bowtie\\
\end{tangle}
\enspace(o) -\enspace
\begin{tangle}
\step\object{A}\step[2]\object{A}\\
\step\XX\\
\dd\step[1]\td \lrcoprod\\
\tu \bowtie\step[2]\id\\
\end{tangle}
\enspace(p) -\enspace
\begin{tangle}
\step\object{A}\step[2]\object{A}\\
\step\XX\\
\dd\step[1]\td \lrcoprod\\
\XX\step[2]\id\\
\id\step[2]\tu \bowtie\\
\end{tangle}\enspace(q) \ \ \ .
\]
By definition, $(n)=(1)-(5)-(12)-(14); (o)=(3)+(10)+(8)+(16);
-(p)=-(2)+(7)+(11)-(13); -(q)=-(4)+(9)-(6)-(15)$. Therefore (LB)
holds on  $A \otimes A$. Similarly, (LB) holds on  $H \otimes H$.
Consequently,  (LB) holds on  $(A \lrbiprod H) \otimes (A \lrbiprod
H)$. $\Box$

In the sequel, we have
\begin {Corollary} \label {cor2.1}
\ (\cite [Proposition 8.3.4] {Ma95}) Assume that $A$ and $H$ are
braided Lie bialgebras;   $(A, H)$ is a matched pair of braided Lie
algebras;  $A$ is a left $H$-module braided Lie coalgebra;
 $H$ is a right $A$-module braided Lie coalgebra. If $(\id_{H}\ot
\alpha)(\delta_{H}\ot \id_{A})+(\beta\ot\id_{A} )(
\id_{H}\ot\delta_{A})=0$ holds, then $A\bowtie H$ becomes a braided
Lie bialgebra.
\end {Corollary}

\begin {Corollary} \label {cor2.2}
 Assume that  $A$ and $H$ are braided Lie bialgebras; $(A, H)$ is a matched pair of  braided Lie coalgebra;
  $A$ is a left $H$-comodule braided Lie algebra;
$H$ is a right $A$-comodule braided Lie algebra. If $([,]_{H}\ot
\id_{A})(\id_{H}\ot \phi)+(\id_{H}\ot [, ]_{A})(\psi \ot \id_{A})=0$
holds, then $A\lrcoprod H$ becomes a braided Lie bialgebra.
\end {Corollary}

\begin {Corollary} \label {cor2.3}
\quad  Assume that $A$ and $H$ are braided Lie bialgebras; $A$ is a
left $H$-comodule braided Lie coalgebra; $H$ is a right $A$-module
braided Lie algebra. If  (B2) and (B3) hold, and $(\beta\ot \id)(
\id \ot \delta_{A})+( [,]_{H}\ot\id)( \id \ot\phi)=0$ holds,
 then $A\lcprp
H$ becomes a braided Lie bialgebra. \end {Corollary}

\begin {Corollary} \label {cor2.4}\quad (\cite[Proposition 8.3.5]{Ma95})
Assume that $A$ and $H$ are braided Lie bialgebras; $A$ is a left
$H$-module braided Lie algebra; $H$ is a right $A$-comodule braided
Lie coalgebra. If  (B1) and (B4) hold,  and $(\id\ot
\alpha)(\delta_{H}\ot \id)+(\id\ot [,]_{A})(\psi \ot \id)=0$ holds,
then $A\lprcp H$ becomes a braided Lie bialgebra.
\end {Corollary}

\begin {Corollary} \label {cor2.5}
\ (\cite [Theorem 3.7]{Ma00})\quad  If $H$ is a braided Lie
bialgebra in $\cal C$ and $A$ is a left $H$-infinitesimal braided
Lie bialgebra, then bisum  $A\lbiprod H$ is a braided Lie bialgebra
in  $\cal C$. \end {Corollary}

\begin {Corollary} \label {cor2.6}
\quad If $A$ is a braided Lie bialgebra in $\cal C$ and
 $H$ is a right $A$-infinitesimal braided Lie bialgebra, then bisum $A\rbiprod
H$ is a braided Lie bialgebra in $\cal C$. \end {Corollary}

\begin {Corollary} \label {cor2.7}\quad
(i) Assume that $A$ and  $H$ are  braided Lie bialgebras in $\cal
C$. If  $A$ is an $H$-module braided Lie algebra and an  $H$-module
braided Lie coalgebra,  and  $(\id\ot \alpha)(\delta_{H}\ot \id) =0$
holds, then $A\lprod H$ is a braided Lie bialgebra in $\cal C$.

(ii) Assume that $A$ and  $H$ are braided Lie bialgebras in $\cal
C$. If  $A$ is an $H$-comodule braided Lie algebra and $H$- comodule
braided Lie coalgebra, and  $( [,]\ot\id)( \id\ot \phi)=0$ holds,
then $A\lcoprod H$ is a braided Lie bialgebra in $\cal C$.
\end {Corollary}

Obviously, for any two braided Lie bialgebra $A$ and $H$, if set
$\alpha =0$, $ \beta=0 $, $\phi =0$ and $\psi =0$, then  $A ^ \phi
_\alpha \bowtie ^\psi _\beta H$ becomes a braided Lie bialgebras. In
this case, denote  $A ^ \phi _\alpha \bowtie ^\psi _\beta H$ by $A
\oplus H$. Therefore, the direct sum $A \oplus H$ of  $A$ and $H$
 is a braided Lie bialgebra. Using this, we can construct many examples of braided Lie bialgebras.

\begin {Example} \label {1.14} Under Example  \ref {4.8}, assume
that  $U$ has a left duality  $U^*$ in ${\cal C}$. Set $X = U
\otimes U^* $ and  $Y = U ^*\otimes U$. Define the bracket operation
and co-bracket operation of  $X \oplus Y$ as follows: $[\ , \ ]_X =
m_X - m _XC_{X, X}$, $\delta _X =0$, $[\ ,  \ ]_Y = 0$, $\delta _X
=\Delta _Y - C _{Y,Y} \Delta _Y$. It is clear that  $X \oplus Y$ is
a braided Lie bialgebra under the following operations: $[\ , \ ]_{X
\oplus Y} = [ \ , \ ] _X, $ $\Delta _{X \oplus Y} = \delta _Y$.
\end {Example}

 \section{The construction of quantum doubles}\label {3}

In this section we construct the quantum double of braided Lie
algebras.

Assume that $H$ is a braided Lie bialgebra in $\cal C$ and $I$ is
the unit object in $\cal C$. If there exists a morphism $R: I\to
H\ot H$ satisfying the following (COB), then $H$ is called a
coboundary braided Lie bialgebra.
\[ (COB): \ \ \
\begin{tangle}
\step\object{H}\\
\cd*\\
\object{H}\step[2]\object{H}\\
\end{tangle}
\enspace =\enspace
\begin{tangle}
\object{H}\\
\id\step[2]\ro R \\
\cu*\step[2]\id\\
\step\object{H}\step[3]\object{H}\\
\end{tangle}
\enspace+\enspace
\begin{tangle}
\object{H}\\
\id\step[2]\ro R \\
\XX\step[2]\id\\
\id\step[2]\cu*\\
\object{H}\step[3]\object{H}\\
\end{tangle} \ \ \ .
\]
If $(H, R)$ is a coboundary braided Lie bialgebra and $R$ satisfies
the following (CYBE), then $(H, R)$ is called quasitriangular
braided Lie bialgebra:
\[
(CYBE): \ \ \ \begin{tangle}
\ro R\step[2]\ro R\\
\id\step[2]\XX\step[2]\id\\
\cu*\step[2]\id\step[2]\id\\
\step\object{H}\step[3]\object{H}\step[2]\object{H}\\
\end{tangle}
\enspace+\enspace
\begin{tangle}
\ro R\step[2]\ro R\\
\id\step[2]\id\step[2]\id\step[2]\id\\
\id\step[2]\cu*\step[2]\id\\
\object{H}\step[3]\object{H}\step[3]\object{H}\\
\end{tangle}
\enspace+\enspace
\begin{tangle}
\ro R\step[2]\ro R\\
\id\step[2]\XX\step[2]\id\\
\id\step[2]\id\step[2]\cu*\\
\object{H}\step[2]\object{H}\step[3]\object{H}\\
\end{tangle} \enspace=0 \ \ \ .
\]
The equation above  is called the classical Yang-Baxter equation. If
 (COB) holds, then  (CYBE) is equivalent to the following (I) and (II), respectively.
\[
\begin{tangle}
\step\ro R\\
\cd*\step\id\\
\object{H}\step[2]\object{H}\step\object{H}\\
\end{tangle}
\enspace=\enspace
\begin{tangle}
\ro R\step[2]\ro R\\
\id\step[2]\XX\step[2]\id\\
\id\step[2]\id\step[2]\cu*\\
\object{H}\step[2]\object{H}\step[3]\object{H}\\
\end{tangle}\enspace\text{(I)} ,
\step[5]
\begin{tangle}
\ro R\\
\id\step\cd*\\
\object{H}\step\object{H}\step[2]\object{H}\\
\end{tangle} \enspace=\enspace
\begin{tangle}
\Ro R\\
\id\step[2]\ro R\d\\
\cu*\step[2]\id\step\id\\
\step\object{H}\step[3]\object{H}\step\object{H}\\
\end{tangle}\enspace\text{(II)} \ \ .
\]
In fact, by  the braided anti- symmetry of $H$ and (COB),
\[
\begin{tangle}
\step\ro R\\
\cd*\step\id\\
\object{H}\step[2]\object{H}\step\object{H}\\
\end{tangle}
\enspace=\enspace
\begin{tangle}
\Ro R\\
\id\step[2]\ro R\d\\
\cu*\step[2]\id\step\id\\
\step\object{H}\step[3]\object{H}\step\object{H}\\
\end{tangle}
\enspace+\enspace
\begin{tangle}
\Ro R\\
\id\step\ro R\step\id\\
\X\step[2]\id\step\id\\
\id\step\cu*\step\id\\
\object{H}\step[2]\object{H}\step[2]\object{H}\\
\end{tangle}
\enspace=-\enspace
\begin{tangle}
\Ro R\\
\id\step[2]\ro R\d\\
\XX\step[2]\id\step\id\\
\cu*\step[2]\id\step\id\\
\step\object{H}\step[3]\object{H}\step\object{H}\\
\end{tangle}
\enspace-\enspace
\begin{tangle}
\Ro R\\
\id\step\ro R\step\id\\
\X\step[2]\id\step\id\\
\id\step\XX\step\id\\
\id\step\cu*\step\id\\
\object{H}\step[2]\object{H}\step[2]\object{H}\\
\end{tangle}
\]
\[
\enspace=-\enspace
\begin{tangle}
\ro R\step[2]\ro R\\
\id\step[2]\XX\step[2]\id\\
\cu*\step[2]\id\step[2]\id\\
\step\object{H}\step[3]\object{H}\step[2]\object{H}\\
\end{tangle}
\enspace-\enspace
\begin{tangle}
\ro R\step[2]\ro R\\
\id\step[2]\id\step[2]\id\step[2]\id\\
\id\step[2]\cu*\step[2]\id\\
\object{H}\step[3]\object{H}\step[3]\object{H}\\
\end{tangle}\ \ \ \ .
\]
Thus (CYBE) and (I) are equivalent. Similarly, (CYBE) and (II) are
equivalent.

Assume that $H$ is a braided Lie bialgabra in $\cal C$ with a left
duality  $H^{*}$. Define the bracket operation and co-bracket
operation in $H^{*}$ as follows:

\[
\begin{tangle}
\object{H^{*}}\step[3]\object{H^{*}}\\
\cu* \step[2]\\
\end{tangle}
\enspace=\enspace
\begin{tangle}

\object{H^{*}}\step[2]\object{H^{*}}\\
\id\step[2]\id\step[3]\coev \\

\id\step[2]\id\step[2]\cd* \step[1]\id\\
\id\step[2]\XX\step[2]\id\step[1]\id\\
\ev\step[2]\ev\step[1]\id\\
 \step[7]\object {H^*}

\end{tangle}
\step[4], \ \ \
\begin{tangle}
\step\object{H^{*}}\\
\cd* \\
\object{H^{*}}\step[2]\object{H^{*}}\\

\end{tangle}
\enspace=\enspace
\begin{tangle}
\object{H^{*}}\\

\id\step[2] \coev  \step[2] \coev \step[2] \\

\id\step[2] \id   \step[2] \XX \step[2]\id\\

\id\step[2] \cu* \step[2] \id \step[2]\id\\

\Ev\step[5] \id \step[2]\id\\
\step[6] \object{H^{*}}\step[2]\object{H^{*}}\\
\end{tangle} \ \ \ .
\]
It is clear that $ H^{*}$ is a braided Lie bialgebra in $\cal C$(
see \cite [Proposition 8.1.2]{Ma95}). Obviously, $H^{*\text{op}}$ is
a braided Lie bialgebra, where
    the co-bracket operations  of $H^{*\text{op}}$ and $H^{*}$ are the
    same, but their bracket operation  have different sign.
Similarly,  $H^{*\text{cop}}$ is also a braided Lie bialgebra.

\begin {Theorem} \label {2} Assume that  $H$  is a braided Lie bialgebra with a left duality  $H^*$.

(i) (\cite [Pro. 8.2.1] {Ma95})  Let  $A = H^{*\text{op}}$,
\[
\begin{tangle}
\object{H}\step[2.5]\object{H^{* \text{op}}}\\
\Ru\\
\object{H}
\end{tangle}
 \enspace=\enspace
\begin{tangle}
\step\object{H}\step[3.5]\object{H^{*\text{op}}}\\
\cd*\step[2]\id\\
\id\step[2]\XX\\
\id\step[2]\ev\\
\object{H}
\end{tangle}
\step[4], \ \
\begin{tangle}
\object{H}\step[2.5]\object{H^{* \text{op}}}\\
\Lu\\
\step[2.5]\object{H^{* \text{op}}}\\
\end{tangle}
 \step=\enspace
\begin{tangle}
\object{H}\step[4]\object{H^{*\text{op}}}\\
\id\step[2]\cd*\\
\XX\step[2]\id\\
\id\step[2]\XX\\
\id\step[2]\ev\\
\object{H^{* \text{op}}}\\
\end{tangle} \ \ \ , \ \ \ R = C_{H, H^*}coev_H =  \ \ \ \  \begin{tangle}
\step[2]\coev\\
\step[2]\XX\\
\object{A \bowtie H}\step[5]\object{ A \bowtie H}
\end{tangle} \ \ \ .
\]
Then  $(A \bowtie H, R)$ is a quasitriangular braided Lie bialgebra.

(ii) Let  $A = H^{*\text{cop}}$,
\[
\begin{tangle}
\object{H}\step[2.5]\object{H^{* \text{cop}}}\\
\Ru\\
\object{H}
\end{tangle}
 \enspace=\enspace
\begin{tangle}

\step\object{H}\step[3.5]\object{H^{*\text{cop}}}\\
\cd*\step[2]\id\\

\id\step[2]\XX\\

\XX \step[2]\id\\

\ev\step[2]\id\\
\step[4]\object{H}
\end{tangle}
\step[4], \ \
\begin{tangle}
\object{H}\step[2.5]\object{H^{* \text{cop}}}\\

\Lu\\
\step[2.5]\object{H^{* \text{cop}}}\\
\end{tangle}
 \step=\enspace
\begin{tangle}
\object{H}\step[4]\object{H^{*\text{cop}}}\\

\id\step[2]\cd*\\
\XX\step[2]\id\\

\ev\step[2]\id\\
\step[4]\object{H^{* \text{cop}}}\\
\end{tangle} \ \ \ , \ \ \ R =  coev_H = \ \ \ \  \begin{tangle}
\step[2]\coev \\
\object{A \bowtie H}\step[5]\object{ A \bowtie H}
\end{tangle} \ \ \ \ \ \ \ .
\]
Then  $(A \bowtie H, R)$ is a quasitriangular braided Lie bialgebra.
$H^{*\text{cop}}\bowtie H$ is called the quantum double\footnote
{Majid \cite [Pro.8.2.1]{Ma95} called
 $H^{*\text {op}} \bowtie H$ the quantum double. Considering the quantum double of Hopf algebras, we do not use this
 definition in this paper} of $H$,
written as $D(H)$ .
\end  {Theorem}

{\bf Proof.}  (i) \text{the left hand side of (M2)} =\[\enspace
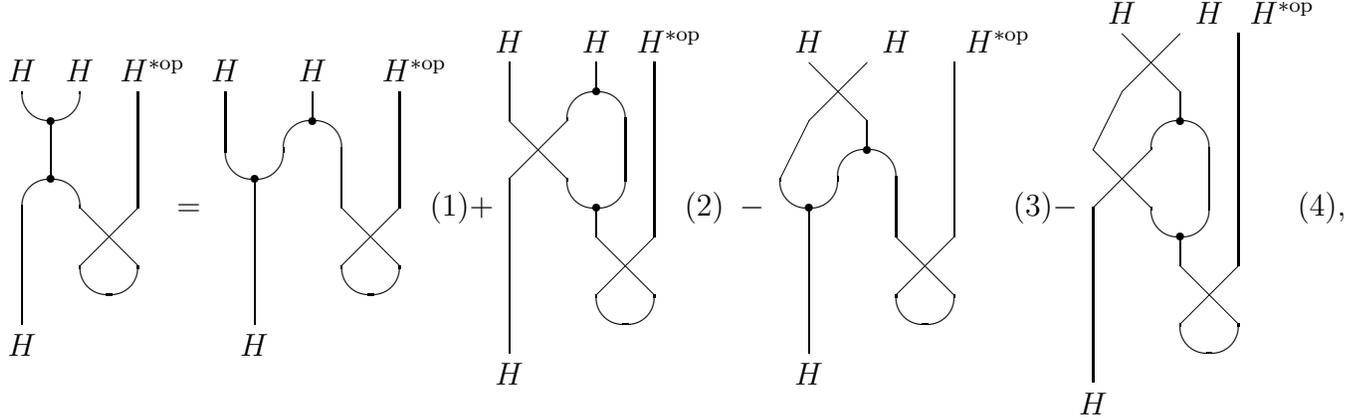
\begin{tangle}
\object{H}\step[2]\object{H}\step[2.5]\object{H^{*\text{op}}}\\
\cu*\step[2]\id\\
\cd*\step[2]\id\\
\id\step[2]\XX\\
\id\step[2]\ev\\
\object{H}
\end{tangle}
\enspace=\enspace
\begin{tangle}
\object{H}\step[3]\object{H}\step[3.5]\object{H^{*\text{op}}}\\
\id\step[2]\cd*\step[2]\id\\
\cu*\step[2]\id\step[2]\id\\
\step\id\step[3]\XX\\
\step\id\step[3]\ev\\
\step\object{H}
\end{tangle}
\enspace(1)+\enspace
\begin{tangle}
\object{H}\step[3]\object{H}\step[2.5]\object{H^{*\text{op}}}\\
\id\step[2]\cd*\step\id\\
\XX\step[2]\id\step\id\\
\id\step[2]\cu*\step\id\\
\id\step[3]\XX\\
\id\step[3]\ev\\
\object{H}
\end{tangle}\enspace(2)
\enspace-\enspace
\begin{tangle}
\step \object{H}\step[3]\object{H}\step[3.5]\object{H^{*\text{op}}}\\
\step [1]\XX\step[3]\id\\

\ne1\step[1]\cd*\step[2]\id\\
\cu*\step[2]\id\step[2]\id\\
\step\id\step[3]\XX\\
\step\id\step[3]\ev\\
\step\object{H}
\end{tangle}
\enspace(3)-\enspace
\begin{tangle}
\step [1]\object{H}\step[3]\object{H}\step[2.5]\object{H^{*\text{op}}}\\
\step [1]\XX\step[2]\id\\
\ne1\step[1]\cd*\step\id\\
\XX\step[2]\id\step\id\\
\id\step[2]\cu*\step\id\\
\id\step[3]\XX\\
\id\step[3]\ev\\
\object{H}
\end{tangle}\enspace(4), \ \
\]
\text{ the right hand side of (M2)}=\[\enspace
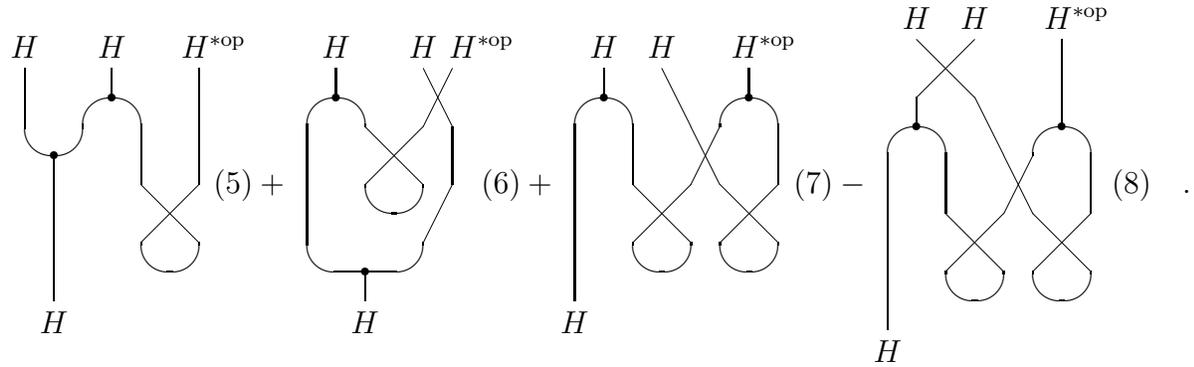
\begin{tangle}
\object{H}\step[3]\object{H}\step[3.5]\object{H^{*\text{op}}}\\

\id\step[2]\cd*\step[2]\id\\
\cu*\step[2]\id\step[2]\id\\
\step\id\step[3]\XX\\
\step\id\step[3]\ev\\
\step\object{H}
\end{tangle}
(5)+\enspace
\begin{tangle}
\step\object{H}\step[3]\object{H}\step[2]\object{H^{*\text{op}}}\\
\cd*\step[2]\X\\
\id\step[2]\XX\step\id\\

\id\step[2]\ev \ne1\\

\Cu*\\
\step[2]\object{H}
\end{tangle}
(6)+\enspace
\begin{tangle}
\step\object{H}\step[2]\object{H}\step[3.5]\object{H^{*\text{op}}}\\
\cd*\step\d\step\cd*\\
\id\step[2]\id\step[2]\X\step[2]\id\\
\id\step[2]\XX\step\XX\\
\id\step[2]\ev\step\ev\\
\object{H}
\end{tangle}
\enspace(7)-\enspace
\begin{tangle}
\step\object{H}\step[2]\object{H}\step[3.5]\object{H^{*\text{op}}}\\
\step\XX\step[3]\id\\
\cd*\step\d\step\cd*\\
\id\step[2]\id\step[2]\X\step[2]\id\\
\id\step[2]\XX\step\XX\\
\id\step[2]\ev\step\ev\\
\object{H}
\end{tangle}\ \  (8) \ \ \ .
\]
It is clear that  (1)= (5), (2)= $-$(8), $-$(3) = (6), (7) =$ -$
(4).  Therefore  (M2) holds. Similarly, we can check that conditions
in Corollary \ref {cor2.1} hold, i.e. $H^{*\text{op}}\bowtie H$ is a
braided Lie bialgebra. See
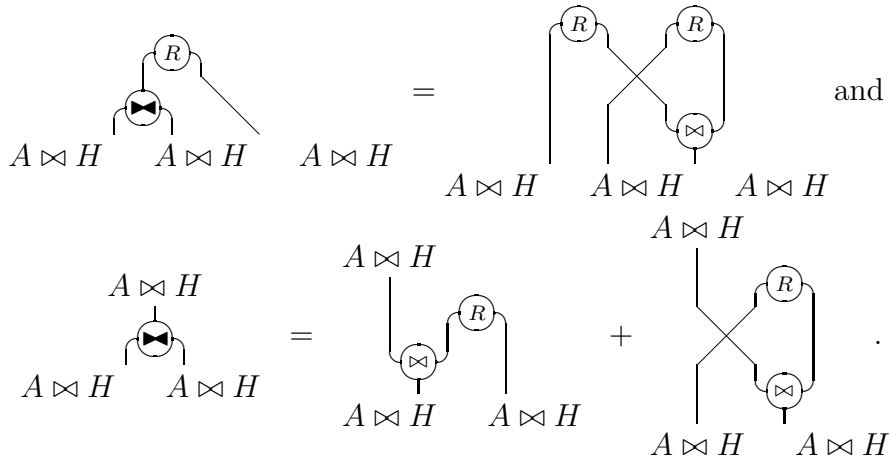
\[
\begin{tangle}
\step[3]\ro R\\
\step[2]\td \lrcoprod \step\nw2\\
\object{A \bowtie H}\step[5]\object{A \bowtie H}\step [5]\object{A \bowtie H}\\
\end{tangle}
\enspace \ \ \ \ = \ \ \ \enspace
\begin{tangle}
\step[2]\ro R\step[2]\ro R\\
\step[2]\id\step[2]\XX\step[2]\id\\
\step[2]\id\step[2]\id\step[2]\tu {\bowtie }\\
\object{A \bowtie H}\step[5]\object{A \bowtie H}\step [5]\object{A \bowtie H}\\
\end{tangle} \ \  \ \ \mbox { and  } \]
\[
\begin{tangle}
\step[2]\step\object{A \bowtie H}\\
\step[2]\td \lrcoprod\\
\object{A \bowtie H}\step[5]\object{A \bowtie H}\\
\end{tangle}
\enspace  \ \ \ \ \ = \ \ \ \ \ \enspace
\begin{tangle}
\object{A \bowtie H}\\
\id\step[2]\ro R \\
\tu {\bowtie }\step[2]\id\\
\object{A \bowtie H}\step[5]\object{A \bowtie H}\\
\end{tangle}
\enspace \ \ \ \ \ + \ \ \ \ \enspace
\begin{tangle}
\object{A \bowtie H}\\
\id\step[2]\ro R \\
\XX\step[2]\id\\
\id\step[2]\tu {\bowtie }\\
\object{A \bowtie H}\step[5]\object{A \bowtie H}\\
\end{tangle} \ \ \ .
 \] Therefore, we complete the proof of the first part. The second part can be proved
 similarly.
 $\Box$

\begin {Corollary} \label {3.4}
 Assume that $(H, R)$ is a quasitriangular braided Lie bialgebra and  $(A, \alpha )$ is an
 $H$-module braided Lie algebra and $H$-module braided Lie coagebra. If we define
 \[ \phi =
\begin{tangle}
\step[2]\object{A}\\
\Ld\\
\object{H}\step[2]\object{A}\\
\end{tangle}
\enspace=\enspace
\begin{tangle}
\step[4]\object{A}\\
\ro R\step[2]\id\\
\XX\step[2]\id\\
\id\step[2]\Lu\\
\object{H}\step[4]\object{A}\\
\end{tangle} \ \ ,
\] and  (SLB) holds on  $A \otimes A$, then  $A\lbiprod H$ becomes a braided Lie bialgebra.\end {Corollary}
{\bf Proof. } We first check that $(A, \phi)$ is an $H$-comodule.
Indeed, by the definition above and the quasitriangular conditions,
\[
\begin{tangle}
\step[4]\object{A}\\
\step\ld[3]\\
\cd*\step[2]\id\\
\object{H}\step[2]\object{H}\step[2]\object{A}
\end{tangle}=\enspace
\begin{tangle}
\step[5]\object{A}\\
\step\ro R\step[2]\id\\
\step\XX\step[2]\id\\
\cd*\step\Lu\\
\object{H}\step[2]\object{H}\step[3]\object{A}\\
\end{tangle}
\enspace=\enspace
\begin{tangle}
\step[6]\object{A}\\
\Ro R\step[2]\id\\
\id\step[2]\ro R\d\step\id\\
\cu*\step\dd\step\id\step\id\\
\step\XX\step[2]\id\step\id\\
\step\id\step[2]\XX\step\id\\
\step\id\step[2]\id\step[2]\lu\\
\step\object{H}\step[2]\object{H}\step[3]\object{A}\\
\end{tangle}
\enspace=\enspace
\begin{tangle}
\step[7]\object{A}\\
\step\Ro R\step[2]\id\\
\dd\step\ro R\step\d\step\id\\
\id\step[2]\XX\step[2]\id\step\id\\
\XX\step[2]\XX\step\id\\
\id\step[2]\XX\step[2]\lu\\
\id\step[2]\id\step[2]\nw2\step[2]\id\\
\id\step[2]\id\step[2]\step[2]\lu\\
\object{H}\step[2]\object{H}\step[5]\object{A}\\
\end{tangle}
\enspace-\enspace
\begin{tangle}
\step[7]\object{A}\\
\step\Ro R\step[2]\id\\
\dd\step\ro R\step\d\step\id\\
\XX\step[2]\id\step[2]\id\step\id\\
\id\step[2]\XX\step[2]\id\step\id\\
\XX\step[2]\XX\step\id\\
\id\step[2]\XX\step[2]\lu\\
\id\step[2]\id\step[2]\nw2\step[2]\id\\
\id\step[2]\id\step[2]\step[2]\lu\\
\object{H}\step[2]\object{H}\step[5]\object{A}\\
\end{tangle}
\]
\[
=\enspace
\begin{tangle}
\step[4]\object{A}\\
\step\ro R\step\id\\
\step\XX\step\id\\
\dd\ro R \lu\\
\id\step\XX\step\id\\
\id\step\id\step[2]\lu\\
\object{H}\step\object{H}\step[3]\object{A}\\
\end{tangle}
\enspace-\enspace
\begin{tangle}
\step[4]\object{A}\\
\step\ro R\step\id\\
\step\XX\step\id\\
\dd\ro R \lu\\
\id\step\XX\step\id\\
\X\step[2]\lu\\
\object{H}\step\object{H}\step[3]\object{A}\\
\end{tangle}
\enspace =\enspace
\begin{tangle}
\step[4]\object{A}\\
\ld[4]\\
\id\step[2]\Ld\\
\object{H}\step[2]\object{H}\step[2]\object{A}\\
\end{tangle}
\enspace-\enspace
\begin{tangle}
\step[4]\object{A}\\
\ld[4]\\
\id\step\step\Ld\\
\XX\step[2]\id\\
\object{H}\step[2]\object{H}\step[2]\object{A}
\end{tangle} \ \ .
\]
Next we check the Yetter-Drinfeld module condition (YD):
\begin{equation*}
\begin{tangle}
\object{H}\step[2]\object{A}\\
\Lu\\
\Ld \\
\object{H}\step[2]\object{A}
\end{tangle}
\enspace=\enspace
\begin{tangle}
\object{H}\step[3]\object{A}\\
\id\step[2]\ld\\
\cu*\step\id\\
\step\object{H}\step[2]\object{A}
\end{tangle}
\enspace(1)+\enspace
\begin{tangle}
\object{H}\step[3]\object{A}\\
\id\step[2]\ld\\
\XX\step\id\\
\id\step[2]\lu\\
\object{H}\step[3]\object{A}
\end{tangle}
\enspace(2)+\enspace
\begin{tangle}
\step\object{H}\step[2]\object{A}\\
\cd*\step\id\\
\id\step[2]\lu\\
\object{H}\step[3]\object{A}
\end{tangle}\enspace(3).
\end{equation*}
Since $A$ is a left $H$-comodule, we have
\[
\text{the left hand side }=\enspace
\begin{tangle}
\step[2]\object{H}\step[2]\object{A}\\
\ro R\lu[2]\\
\XX\step[2]\id\\
\id\step\step\lu[2]\\
\object{H}\step[4]\object{A}\\
\end{tangle}
\enspace=\enspace
\begin{tangle}
\step[4]\object{H}\step\object{A}\\
\ro R\step[2]\id\step\id\\
\XX\step[2]\id\step\id\\
\id\step[2]\cu*\step\id\\
\id\step[3]\Lu\\
\object{H}\step[5]\object{A}\\
\end{tangle}
\enspace(a)+\enspace
\begin{tangle}
\step[4]\object{H}\step\object{A}\\
\ro R\step[2]\id\step\id\\
\XX\step[2]\id\step\id\\
\id\step[2]\XX\step\id\\
\id\step[2]\nw2\step\lu\\
\id\step[4]\lu\\
\object{H}\step[5]\object{A}\\
\end{tangle}\enspace(b).
\]
For the right hand side,
\[
(3)=-\enspace
\begin{tangle}
\step\object{H}\step[2]\object{A}\\
\cd*\step\id\\
\XX\step\id\\
\id\step[2]\lu\\
\object{H}\step[3]\object{A}
\end{tangle}
\enspace=-
\begin{tangle}
\object{H}\step[5]\object{A}\\
\id\step[2]\ro R\step\id\\
\cu*\step\dd\step\id\\
\step\XX\step[2]\id\\
\step\id\step[2]\Lu\\
\step\object{H}\step[4]\object{A}
\end{tangle} \ \ \ (4)
\enspace-\enspace
\begin{tangle}
\object{H}\step[5]\object{A}\\
\id\step[2]\ro R\step\id\\
\XX\step[2]\id\step\id\\
\d\step\cu*\step\id\\
\step\XX\step[2]\id\\
\step\id\step[2]\Lu\\
\step\object{H}\step[4]\object{A}
\end{tangle} \ \ \ (5) \ ,
\enspace
\]
where $-(4)=(a)$, $-(5) +(1)=0$, $(2)=(b)$.  Therefore $(YD) $
holds. Using Corollary  \ref {cor2.5}, we complete the proof. $\Box$

Dually, we shall give the definition of co-quasitriangular braided
Lie bialgebras. Assume that $B$ is a braided Lie bialgebra. If there
exists a morphism $r: H\ot H\to I$ satisfying the following
condition  (Bo), then $(H, r)$ is called a coboundary Lie bialgebra.
Furthermore, if $r$ satisfies the following condition (CCYBE), then
$(H, r)$ is called a co-quasitriangular braided Lie bialgebra:
\[ (Bo): \ \ \
\begin{tangle}
\object{H}\step[2]\object{H}\\
\cu*\\
\step\object{H}
\end{tangle}
\enspace =\enspace
\begin{tangle}
\step\object{H}\step[3]\object{H}\\
\cd*\step[2]\id\\
\id\step[2]\coro r \\
\object{H}
\end{tangle}
\enspace+\enspace
\begin{tangle}
\object{H}\step[3]\object{H}\\
\id\step[2]\cd*\\
\x\step[2]\id\\
\id\step[2]\coro r \\
\object{H}
\end{tangle} \ \ \ . \]
\[ (CCYBE): \ \ \
\begin{tangle}
\step\object{H}\step[3]\object{H}\step[2]\object{H}\\
\cd*\step[2]\id\step[2]\id\\
\id\step[2]\x\step[2]\id\\
\coro r\step[2]\coro r\\
\end{tangle}
\enspace+\enspace
\begin{tangle}
\object{H}\step[3]\object{H}\step[3]\object{H}\\
\id\step[2]\cd*\step[2]\id\\
\id\step[2]\id\step[2]\id\step[2]\id\\
\coro r\step[2]\coro r\\
\end{tangle}
\enspace+\enspace
\begin{tangle}\object{H}\step[2]\object{H}\step[3]\object{H}\\
\id\step[2]\id\step[2]\cd*\\
\id\step[2]\x\step[2]\id\\
\coro r\step[2]\coro r\\
\end{tangle} \enspace= 0 \ \ \ .\]
The equation above is called  the classical Co-Yang-Baxter equation.


If (Bo) holds, then  (CCYBE)  is equivalent to anyone of the
following (I) and  (II):
\[
\begin{tangle}
\object{H}\step[2]\object{H}\step\object{H}\\
\cu*\step\id\\\step\coro r\end{tangle} \enspace=\enspace
\begin{tangle}\object{H}\step[2]\object{H}\step[3]\object{H}\\
\id\step[2]\id\step[2]\cd*\\
\id\step[2]\x\step[2]\id\\
\coro r\step[2]\coro r\end{tangle}  \ \ \ (I) \ ,  \ \ \ \
\begin{tangle}
\object{H}\step\object{H}\step[2]\object{H}\\
\id\step\cu*\\
\coro r\end{tangle} \enspace=\enspace
\begin{tangle}
\step\object{H}\step[3]\object{H}\step\object{H}\\
\cd*\step[2]\id\step\id\\
\id\step[2]\coro r\dd\\
\coRo r\\
\end{tangle} \ \ \ (II) \ \ \ .
\]

\begin {Corollary} \label {3.5}
Assume  that $(H, r)$ is a co-quasitriangular braided Lie bialgebra;
$(A, \phi)$ is an
 $H$-comodule braided Lie algebra and an $H$-comodule braided Lie coalgebra. If we define
 \[ \alpha =\ \
\begin{tangle}
\object{H}\step[2]\object{A}\\
\Lu\\
\step[2]\object{A}\\
\end{tangle}
\enspace=\enspace
\begin{tangle}
\object{H}\step[4]\object{A}\\
\id\step[2]\Ld\\
\XX\step[2]\id\\
\coro r\step[2]\id\\
\step[4]\object{A}\\
\end{tangle} \ \ \ ,
\]
 and  (SLB) holds on $A \otimes A$, then  $A\lbiprod H$ becomes a braided Lie bialgebra.
\end {Corollary}

Corollary \ref {3.4} and Corollary \ref {3.5} are called the
bosonisation theorems.

We easily obtain: if $H$ has a left duality $H^*$, then   $(H, R)$
is a quasitriangular braided Lie bialgebra if and only if
 $(H ^*,
R^*)$ is a co-quasitriangular braided Lie bialgebra. Here  $R^* $
can be denoted by the following diagram:
\[ R^* = \ \ \ \begin{tangle}
\object{H^* }\step[2]\object{H^*}\\
\id\step[2]\id \step [2] \ro R\\

\id \step [2]\XX\step[2]\id\\
\ev \step[2]\ev \step [2]\\
\end{tangle} \ \ \ .
\]

Otherwise, if  $U$ and  $V$ have the left duality  $U^*$ and $V^*$,
then  $U \oplus V$ has a left duality  $U^*\oplus V^*$. In fact, let
$D $ denote  $ U \oplus V$. Define  $\widetilde{ d_U }= d_U (\pi
_{U^*} \otimes \pi _U) $, $\widetilde { b_U } $ $= (\iota _U
\otimes\iota _{U^*})b_U$,  $d_D = \widetilde{d_U }+ \widetilde{d
_V}$, $b_D = \widetilde {b_U}+ \widetilde {b_V}$, where
 $\pi _{U}$
and  $\pi _{U^*}$ are the canonical projections from $U \oplus U^*$
to  $U$  and to  $U^*$, respectively;  $ \iota_ U$ and $\iota_{U^*}$
 are canonical injection from $U$ to  $U \oplus U^*$  and from  $U^*$ to  $U \oplus U^*$, respectively.
It is clear that  $d_D$ and $b_D$ is evaluation and coevaluation,
respectively. Therefore,  in Example \ref {1.14}, $X \oplus Y$ has a
left duality, so $X \oplus Y$ has the quantum double,

\begin {Example} \label {4.9} (\cite[Lemma 3.4]{AS98b}) Assume that $\Gamma$ is a finite commutative
group; $g_i \in \Gamma,$ $\chi _i\in \hat \Gamma$, where $\hat
\Gamma$ is the character group of $\Gamma$;
 $\chi _i (g_j)\chi_j
(g_i) =1$, $1< N_i $, where  $N_i$ is the order of  $\chi _i (g_i)$,
$1 \le i, j \le \theta$. Let ${\cal R} ( g_i, \chi_i; 1 \le i \le
\theta)$ be the algebra generated by  the set $ \{ x_i \mid 1 \le i
\le \theta \}$ with defining relations
\begin {eqnarray} \label {qlse1} x_l ^{N_l}=0,\  x_{i}x_{j} = \chi_{j}(g_{i})
x_{j}x_{i}  \ \ \ \hbox { where }  1\le i, j,l \le \theta.
\end {eqnarray}  Define the coalgebra operations and  $k\Gamma$-(co)module in ${\cal R} ( g_i, \chi_i; 1 \le i \le
\theta)$ as follows;
$$\Delta x_i = x_i \otimes 1 + 1 \otimes x_i, \ \ \epsilon (x_i)
=0,$$
 $$\delta ^-(x_{i}) = g_{i}\otimes x_{i}, \qquad h \cdot x_{i} =
\chi_{i}(h)x_{i}. $$  It is clear that  ${\cal R} ( g_i, \chi_i; 1
\le i \le \theta)$ is a quantum commutative braided Hopf algebra in
 $^{k\Gamma}_{k\Gamma} {\cal YD} $, called  a quantum linear
  space in
  $^{k\Gamma}_{k\Gamma} {\cal YD} $. Let   $ {\cal C} := ^{k\Gamma}_{k\Gamma} {\cal
YD}$ and   $U := {\cal R} ( g_i, \chi_i; 1 \le i \le \theta)$ as
braided Hopf algebras; $A = H =: {\cal R} ( g_i, \chi_i; 1 \le i \le
\theta)$ as vector spaces.
  Define  $\delta _A =0$, $[\ , \ ]_H = 0, $
$[\ \ ]_A = m_U - m _UC_{U, U},$  \ \ $ \delta _H= \Delta_U - C_{U,
U}\Delta _U .$  It is clear that  $A$ and  $H$  are braided
bialgebras in  $ {\cal C}$ (see \cite [Theorem 1.3]{ZZ04}). So is $A
\oplus H$. Since $A$ and  $H$ are finite dimensional, they have left
dualities  in $^{k\Gamma}_{k\Gamma}{\cal YD}$ (see \cite
[Proposition 1.0.17]{Zh99} ), so  $A \oplus H$ has a left duality in
$^{k\Gamma}_{k\Gamma}{\cal YD}$, i.e. $A \oplus H$ has a quantum
double in $^{k\Gamma}_{k\Gamma}{\cal YD}$.

\end {Example}

 \section {
The universal enveloping algebra of double cross sum of Lie
bialgebras }\label {s4}

In this section  we show that the universal enveloping algebra of
double cross sum of Lie bialgebras is isomorphic to the double cross
product of their universal enveloping algebras.

Throughout this section, we work over ordinary vector spaces. That
is, ${\cal C}$ is the category ${\cal V}ect_k $  of ordinary vector
spaces (see \cite {Zh99}). Let  $i_H\ : H \rightarrow U(H) \ \ $
denote the canonical injection. One can see \cite [Definition
IX.2.1] {Ka95} and \cite {ZC01} about the matched pair of
bialgebras.

\begin {Lemma} \label {4.2}
 Assume that  $A$ and $ H$ are two Lie algebras. If  $(U(A), U(H), \alpha ', \beta')$ is a
 matched pair and    $Im (\alpha '(i_H\otimes i_A))$ $
\subseteq A$ and  $Im (\beta '(i_H\otimes i_A))$ $ \subseteq H$,
 then there exist  $\alpha$ and  $\beta $ such that  $(A, H, \alpha, \beta
)$ becomes a matched pair of Lie algebras with  $\alpha '(i_H
\otimes i_A) = i_A \alpha $ and  \ $\beta '(i_H \otimes i_A) = i_H
\beta $.
\end {Lemma}

{\bf Proof.}  Define  $\alpha(x\ot a)=\alpha'(x\ot a)= x\trr a$ and
$ \beta(x\ot a)= \beta'(x\ot a)=x\trl a$ for any  $x\in H$, $a\in
A$.  Since ($U(A)$ and $U(H)$) are a matched pair of bialgebras, we
have $x\trr(ab)=(x\trr a)((1\trl 1)\trr b)+(1\trr a)((x\trl 1)\trr
b)+(1\trr 1)((x\trl a)\trr b)=(x\trr a)b+a(x\trr b)+(x\trl a)\trr b$
and $x\trr(ba)=(x\trr b)a+b(x\trr a)+(x\trl b)\trr a$ for any  $x\in
H$, $a, b\in A$. By subtraction of two equations above,  $x\trr[a,
b]=[x\trr a, b] +[a, x\trr b]+(x\trl a)\trr b-(x\trl b)\trr a$, so
(M1) holds. Similarly, we can obtain (M2): $[x, y]\trl a=[x, y\trl
a]+[x\trl a, y]+x\trl (y\trr a)-y\trl (x\trr a)$, Consequently, $(A,
H)$ is a matched pair of Lie algebras. $\Box$

\begin {Theorem} \label {3} Assume that $A$ and $ H$ are Lie algebras. If
 $(U(A), U(H), \alpha ',
\beta')$  is a matched pair of bialgebras and  $Im (\alpha
'(i_H\otimes i_A)) $ $ \subseteq A$, $Im (\beta '(i_H\otimes i_A)) $
$\subseteq H$,  then there exist a Hopf algebra isomorphism:
$$ U(A \bowtie
H)\cong U(A)\bowtie U(H). $$  \end {Theorem}

{\bf Proof. } Assume that $(a, x)$,  $ (b, y)\in A\bowtie H$ with
$a, b\in A$, $x, y\in H$. By the definition of $A\bowtie H$, $[(a,
x), (b, y)]=([a, b]+x\trr b-y\trr a, [x,y]+x\trl b- y\trl a)$.
Define $f: A\bowtie H\to U(A)\bowtie U(H), f(a, x)=i_{A}(a)\ot
1+1\ot i_{H}(x)= a\ot 1+1\ot x$. Then  $f([(a, x), (b, y)])=[a,
b]\ot 1+x\trr b\ot 1- y\trr a\ot 1+
 1\ot [x, y]+1\ot x\trl b- 1\ot y\trl a)$.
By the definition of double cross product ( see \cite{ZC01}), $(a\ot
x)\cdot_{\bowtie}(b\ot y)=(a\bowtie x)(b\bowtie y)=a(x\di\trr
b\di)\ot(x\dii\trl b\dii)y$, where $\cdot_{\bowtie}$ denotes the
multiplication of double cross product. See
\begin{align*}
&[f(a, x), f(b, y)]\\
=&\ (a\ot 1+1\ot x)\cdot_{\bowtie}(b\ot 1+1\ot y)-(b\ot 1+1\ot y)\cdot_{\bowtie}(a\ot 1+1\ot x)\\
=&\ ab\ot 1+a\ot y+b\ot x+x\trr b\ot 1
+1\ot x\trl b+1\ot xy\\
&-(ba\ot 1+b\ot x+a\ot y+y\trr a\ot 1
+ 1\ot y\trl a+ 1\ot yx)\\
=&\ (ab-ba)\ot 1+x\trr b\ot 1-y\trr a\ot 1\\
&+ 1\ot (xy-yx)+1\ot x\trl b- 1\ot y\trl a).
\end{align*}
Thus $f([(a, x), (b, y)])=[f(a, x), f(b, y)]$, which implies that
$f$   is a Lie algebra homomorphism from $A\bowtie H$ to
$U(A)\bowtie U(H)$. Using the universal property of universal
enveloping algebras, we know that $f$ can become
 an algebra homomorphism from $U(A\bowtie H)$ to $U(A)\bowtie
U(H)$.

 $f$ also is a coalgebra homomorphism. Indeed, for any $(a, x)\in A\bowtie
H$, since
\begin{eqnarray*}
\Delta f(a, x)&=&\Delta(a\ot 1+1\ot x)=a\ot 1\ot 1\ot 1\\
&&+1\ot 1\ot a\ot 1
+1\ot x\ot 1\ot 1+1\ot 1\ot 1\ot x \ \ \ \mbox { and  }\\
(f\ot f)\Delta (a, x)&=&(f\ot f)((a, x)\ot 1 +1\ot (a, x))\\
&=&(f\ot f)(a\ot 1+x\ot 1+1\ot a +1\ot x)\\
&=&f(a)\ot f(1)+f(x)\ot f(1)+f(1)\ot f(a) +f(1)\ot f(x)),
\end{eqnarray*}
we have $\Delta f(a, x) = (f\ot f)\Delta (a, x)$. Considering that
$f$ and  $\Delta$ are algebra homomorphisms, we have that $f$ is a
coalgebra homomorphism.

 Using two Lie algebra homomorphisms $\iota_{A}: A\to A\bowtie H$ ,   $a\to (a,
 0)$ and
$\iota_{H}: H\to A\bowtie H$, $x\to (0, x)$, we can obtain two
algebra homomorphisms $j_{1}: U(A)\to U(A\bowtie H)$ and  $j_{2}:
U(H)\to U(A\bowtie H)$, i.e. for any $a\in A$, $x\in H$,
$$j_1(a)=i_{A\bowtie H}(a, 0), \quad j_2(x)=i_{A\bowtie H}(0, x).$$
Define a map $j: U(A)\bowtie U(H)\to U(A\bowtie H)$ by sending
$j(b\ot y)=j_1(b)j_2(y)$

For any  $a\in A, x \in H$, since $jf(a, x)=j(a\ot 1)+j(1\ot x)=(a,
0)+(0, x)=(a, x)$, we know that    $f$ is injective on $A \bowtie
H$. By \cite [Lemma 11.0.1]{Sw69}, $f$  is injective on  $U(A
\bowtie H)$.  Obviously, $f$ is surjective. Consequently $f$ is a
Hopf algebra isomorphism.  $\Box$

\begin {Corollary} \label {3.4} If $A$ and  $H$ are Lie algebras, then
there exists a Hopf algebra isomorphism:
$$ U(A \oplus
H)\cong U(A)\otimes U(H). $$  \end {Corollary}

{\bf Proof.} Let  $\alpha '$ and  $\beta '$ be trivial actions, i.e.
$\alpha ' = (\epsilon _{U(H)} \otimes id _{U(A)})$ and  $\beta ' = (
id _{U(H)}\otimes \epsilon _{U(A)})$. Obviously,  $(U(A), U(H),
\alpha ', \beta')$ is a matched pair of bialgebras with
 $\alpha =0$ and  $\beta =0$. Applying Theorem \ref
{3} we complete the proof. $\Box$

\end{document}